\catcode`\@=11

\font\tenmsa=msam10 \font\sevenmsa=msam7 \font\fivemsa=msam5
\font\tenmsb=msbm10
\font\sevenmsb=msbm7 \font\fivemsb=msbm5 \newfam\msafam
\newfam\msbfam
\textfont\msafam=\tenmsa \scriptfont\msafam=\sevenmsa
\scriptscriptfont\msafam=\fivemsa \textfont\msbfam=\tenmsb
\scriptfont\msbfam=\sevenmsb \scriptscriptfont\msbfam=\fivemsb

\def\hexnumber@#1{\ifnum#1<10 \number#1\else \ifnum#1=10
A\else\ifnum#1=11
 B\else\ifnum#1=12 C\else \ifnum#1=13 D\else\ifnum#1=14
E\else\ifnum#1=15
 F\fi\fi\fi\fi\fi\fi\fi}

\def\msa@{\hexnumber@\msafam} \def\msb@{\hexnumber@\msbfam}
\mathchardef\boxdot="2\msa@00 \mathchardef\boxplus="2\msa@01
\mathchardef\boxtimes="2\msa@02 \mathchardef\square="0\msa@03
\mathchardef\blacksquare="0\msa@04 \mathchardef\centerdot="2\msa@05
\mathchardef\lozenge="0\msa@06 \mathchardef\blacklozenge="0\msa@07
\mathchardef\circlearrowright="3\msa@08
\mathchardef\circlearrowleft="3\msa@09
\mathchardef\rightleftharpoons="3\msa@0A
\mathchardef\leftrightharpoons="3\msa@0B
\mathchardef\boxminus="2\msa@0C
\mathchardef\Vdash="3\msa@0D \mathchardef\Vvdash="3\msa@0E
\mathchardef\vDash="3\msa@0F \mathchardef\twoheadrightarrow="3\msa@10
\mathchardef\twoheadleftarrow="3\msa@11
\mathchardef\leftleftarrows="3\msa@12
\mathchardef\rightrightarrows="3\msa@13
\mathchardef\upuparrows="3\msa@14
\mathchardef\downdownarrows="3\msa@15
\mathchardef\upharpoonright="3\msa@16

\mathchardef\downharpoonright="3\msa@17
\mathchardef\upharpoonleft="3\msa@18
\mathchardef\downharpoonleft="3\msa@19
\mathchardef\rightarrowtail="3\msa@1A
\mathchardef\leftarrowtail="3\msa@1B
\mathchardef\leftrightarrows="3\msa@1C
\mathchardef\rightleftarrows="3\msa@1D
\mathchardef\Lsh="3\msa@1E \mathchardef\Rsh="3\msa@1F
\mathchardef\rightsquigarrow="3\msa@20
\mathchardef\leftrightsquigarrow="3\msa@21
\mathchardef\looparrowleft="3\msa@22
\mathchardef\looparrowright="3\msa@23 \mathchardef\circeq="3\msa@24
\mathchardef\succsim="3\msa@25 \mathchardef\gtrsim="3\msa@26
\mathchardef\gtrapprox="3\msa@27 \mathchardef\multimap="3\msa@28
\mathchardef\therefore="3\msa@29 \mathchardef\because="3\msa@2A
\mathchardef\doteqdot="3\msa@2B 
\mathchardef\traceiangleq="3\msa@2C \mathchardef\precsim="3\msa@2D
\mathchardef\lesssim="3\msa@2E \mathchardef\lessapprox="3\msa@2F
\mathchardef\eqslantless="3\msa@30 \mathchardef\eqslantgtr="3\msa@31
\mathchardef\curlyeqprec="3\msa@32 \mathchardef\curlyeqsucc="3\msa@33
\mathchardef\preccurlyeq="3\msa@34 \mathchardef\leqq="3\msa@35
\mathchardef\leqslant="3\msa@36 \mathchardef\lessgtr="3\msa@37
\mathchardef\backprime="0\msa@38 \mathchardef\risingdotseq="3\msa@3A
\mathchardef\fallingdotseq="3\msa@3B
\mathchardef\succcurlyeq="3\msa@3C
\mathchardef\geqq="3\msa@3D \mathchardef\geqslant="3\msa@3E
\mathchardef\gtrless="3\msa@3F \mathchardef\sqsubset="3\msa@40
\mathchardef\sqsupset="3\msa@41
\mathchardef\trianglerighteq="3\msa@44
\mathchardef\trianglelefteq="3\msa@45 \mathchardef\bigstar="0\msa@46
\mathchardef\between="3\msa@47
\mathchardef\blacktriangledown="0\msa@48
\mathchardef\blacktriangleright="3\msa@49
\mathchardef\blacktriangleleft="3\msa@4A
\mathchardef\blacktriangle="0\msa@4E
\mathchardef\triangledown="0\msa@4F
\mathchardef\eqcirc="3\msa@50 \mathchardef\lesseqgtr="3\msa@51
\mathchardef\gtreqless="3\msa@52 \mathchardef\lesseqqgtr="3\msa@53
\mathchardef\gtreqqless="3\msa@54 \mathchardef\Rrightarrow="3\msa@56
\mathchardef\Lleftarrow="3\msa@57 \mathchardef\veebar="2\msa@59
\mathchardef\barwedge="2\msa@5A \mathchardef\doublebarwedge="2\msa@5B
\mathchardef\angle="0\msa@5C \mathchardef\measuredangle="0\msa@5D
\mathchardef\sphericalangle="0\msa@5E
\mathchardef\varpropto="3\msa@5F
\mathchardef\smallsmile="3\msa@60 \mathchardef\smallfrown="3\msa@61
\mathchardef\Subset="3\msa@62 \mathchardef\Supset="3\msa@63
\mathchardef\Cup="2\msa@64 
\mathchardef\Cap="2\msa@65
 \mathchardef\curlywedge="2\msa@66
\mathchardef\curlyvee="2\msa@67 \mathchardef\leftthreetimes="2\msa@68
\mathchardef\rightthreetimes="2\msa@69
\mathchardef\subseteqq="3\msa@6A
\mathchardef\supseteqq="3\msa@6B \mathchardef\bumpeq="3\msa@6C
\mathchardef\Bumpeq="3\msa@6D \mathchardef\lll="3\msa@6E

\mathchardef\ggg="3\msa@6F 
\mathchardef\circledS="0\msa@73
\mathchardef\pitchfork="3\msa@74 \mathchardef\dotplus="2\msa@75
\mathchardef\backsim="3\msa@76 \mathchardef\backsimeq="3\msa@77
\mathchardef\complement="0\msa@7B \mathchardef\intercal="2\msa@7C
\mathchardef\circledcirc="2\msa@7D \mathchardef\circledast="2\msa@7E
\mathchardef\circleddash="2\msa@7F
\def\ulcorner{\delimiter"4\msa@70\msa@70 }
\def\urcorner{\delimiter"5\msa@71\msa@71 }
\def\llcorner{\delimiter"4\msa@78\msa@78 }
\def\lrcorner{\delimiter"5\msa@79\msa@79 }
\def\yen{\mathhexbox\msa@55 }
\def\checkmark{\mathhexbox\msa@58 } \def\circledR{\mathhexbox\msa@72
}
\def\maltese{\mathhexbox\msa@7A } \mathchardef\lvertneqq="3\msb@00
\mathchardef\gvertneqq="3\msb@01 \mathchardef\nleq="3\msb@02
\mathchardef\ngeq="3\msb@03 \mathchardef\nless="3\msb@04
\mathchardef\ngtr="3\msb@05 \mathchardef\nprec="3\msb@06
\mathchardef\nsucc="3\msb@07 \mathchardef\lneqq="3\msb@08
\mathchardef\gneqq="3\msb@09 \mathchardef\nleqslant="3\msb@0A
\mathchardef\ngeqslant="3\msb@0B \mathchardef\lneq="3\msb@0C
\mathchardef\gneq="3\msb@0D \mathchardef\npreceq="3\msb@0E
\mathchardef\nsucceq="3\msb@0F \mathchardef\precnsim="3\msb@10
\mathchardef\succnsim="3\msb@11 \mathchardef\lnsim="3\msb@12
\mathchardef\gnsim="3\msb@13 \mathchardef\nleqq="3\msb@14
\mathchardef\ngeqq="3\msb@15 \mathchardef\precneqq="3\msb@16
\mathchardef\succneqq="3\msb@17 \mathchardef\precnapprox="3\msb@18
\mathchardef\succnapprox="3\msb@19 \mathchardef\lnapprox="3\msb@1A
\mathchardef\gnapprox="3\msb@1B \mathchardef\nsim="3\msb@1C
\mathchardef\napprox="3\msb@1D
\mathchardef\nsubseteqq="3\msb@22
\mathchardef\nsupseteqq="3\msb@23 \mathchardef\subsetneqq="3\msb@24
\mathchardef\supsetneqq="3\msb@25
\mathchardef\subsetneq="3\msb@28
\mathchardef\supsetneq="3\msb@29 \mathchardef\nsubseteq="3\msb@2A
\mathchardef\nsupseteq="3\msb@2B \mathchardef\nparallel="3\msb@2C
\mathchardef\nmid="3\msb@2D \mathchardef\nshortmid="3\msb@2E
\mathchardef\nshortparallel="3\msb@2F \mathchardef\nvdash="3\msb@30
\mathchardef\nVdash="3\msb@31 \mathchardef\nvDash="3\msb@32
\mathchardef\nVDash="3\msb@33 \mathchardef\ntrianglerighteq="3\msb@34
\mathchardef\ntrianglelefteq="3\msb@35
\mathchardef\ntriangleleft="3\msb@36
\mathchardef\ntriangleright="3\msb@37
\mathchardef\nleftarrow="3\msb@38
\mathchardef\nrightarrow="3\msb@39 \mathchardef\nLeftarrow="3\msb@3A
\mathchardef\nRightarrow="3\msb@3B
\mathchardef\nLeftrightarrow="3\msb@3C
\mathchardef\nleftrightarrow="3\msb@3D
\mathchardef\divideontimes="2\msb@3E
\mathchardef\varnothing="0\msb@3F \mathchardef\nexists="0\msb@40
\mathchardef\mho="0\msb@66 \mathchardef\thorn="0\msb@67
\mathchardef\beth="0\msb@69 \mathchardef\gimel="0\msb@6A
\mathchardef\daleth="0\msb@6B \mathchardef\lessdot="3\msb@6C
\mathchardef\gtrdot="3\msb@6D \mathchardef\ltimes="2\msb@6E
\mathchardef\rtimes="2\msb@6F \mathchardef\shortmid="3\msb@70
\mathchardef\shortparallel="3\msb@71
\mathchardef\smallsetminus="2\msb@72
\mathchardef\thicksim="3\msb@73 \mathchardef\thickapprox="3\msb@74
\mathchardef\approxeq="3\msb@75 \mathchardef\succapprox="3\msb@76
\mathchardef\precapprox="3\msb@77
\mathchardef\curvearrowleft="3\msb@78
\mathchardef\curvearrowright="3\msb@79 \mathchardef\digamma="0\msb@7A
\mathchardef\varkappa="0\msb@7B \mathchardef\hslash="0\msb@7D
\mathchardef\hbar="0\msb@7E \mathchardef\backepsilon="3\msb@7F
\def\Bbb{\ifmmode\let\next\Bbb@\else
\def\next{\errmessage{Use \string\Bbb\space only in math
mode}}\fi\next}
\def\Bbb@#1{{\Bbb@@{#1}}} \def\Bbb@@#1{\fam\msbfam#1}

\catcode`\@=\active

\def\del{\partial}

\def\CF{\hbox{{$\cal F$}}} \def\CG{\hbox{{$\cal G$}}}

\def\CR{\hbox{{$\cal R$}}} 
 
\def\CN{\hbox{{$\cal N$}}}

\def\CQ{\hbox{{$\cal Q$}}}

\def\cg{\hbox{{\sl g}}} 

\def\lform{\hbox{$\sqcup$}\llap{\hbox{$\sqcap$}}}

\def\h{{{1\over2}}}

\def\R{{\Bbb R}}
\def\C{{\Bbb C}}
\def\Z{{\Bbb Z}}
\def\N{{\Bbb N}}

\def\eps{{\epsilon}}

\def\lcross{{>\!\!\!\triangleleft}}

\def\bicross{{\blacktriangleright\!\!\!\triangleleft}}
\def\dcross{{\bowtie}}

\def\rbiprod{{\cdot\kern-.33em\triangleright\!\!\!<}}
\def\lbiprod{{>\!\!\!\triangleleft\kern-.33em\cdot\, }}

\def\tens{\mathop{\otimes}}

\def\la{{\triangleright}}\def\ra{{\triangleleft}}

\def\extd{{{\rm d}}}

\def\isom{{\cong}}

\def\span{{\rm span}}

\def\Ad{{\rm Ad}}

\def\id{{\rm id}}

\def\<{\langle}
\def\>{\rangle}

\def\equad{\kern -1.7em}

\def\eqn#1#2{\begin{equation}#2\label{#1}\end{equation}}

\def\s#1{{}_{\scriptscriptstyle(#1)}}
\def\o{{}_{\scriptscriptstyle(1)}}
\def\t{{}_{\scriptscriptstyle(2)}}
\def\th{{}_{\scriptscriptstyle(3)}}
\def\fo{{}_{\scriptscriptstyle(4)}}
\def\fiv{{}_{\scriptscriptstyle(5)}}
\def\six{{}_{\scriptscriptstyle(6)}}
\def\sev{{}_{\scriptscriptstyle(7)}}

\def\bo{{}^{\bar{\scriptscriptstyle(1)}}}
\def\bt{{}^{\bar{\scriptscriptstyle(2)}}}

\def\und#1{{\underline {#1}}}

\def\Bo{{{}_{\und{\scriptscriptstyle(1)}}}}
\def\Bt{{{}_{\und{\scriptscriptstyle(2)}}}}

\def\text#1{\mbox{\rm #1}}
\def\note#1{}

\def\blacksquare{{\lform}}
\def\frac#1#2{{{#1\over#2}}}

\def\proof{\goodbreak\noindent{\bf Proof\quad}}

\def\endproof{{\ $\lform$}\bigskip }

\def\align#1{\begin{eqnarray*}#1\end{eqnarray*}}

\def\cmath#1{\[\begin{array}{c} #1 \end{array}\]}

\documentstyle[11pt,amscd]{article}
\newtheorem{lemma}{Lemma}[section]
\newtheorem{propos}[lemma]{Proposition}
\newtheorem{example}[lemma]{Example}
\newtheorem{theorem}[lemma]{Theorem}

\newtheorem{corol}[lemma]{Corollary}
\newtheorem{defin}[lemma]{Definition}
\newtheorem{remark}[lemma]{Remark}

\begin{document}\baselineskip 20pt

{\ }\qquad\qquad \hskip 4.3in  DAMTP/97-60
\vspace{.2in}

\begin{center} {\LARGE QUANTUM DIFFERENTIALS AND THE
$q$-MONOPOLE REVISITED\footnote{Research
supported by the EPSRC grant GR/K02244}}
\\ \baselineskip 13pt{\ }
{\ }\\
Tomasz Brzezi\'nski\\
+\\
Shahn Majid\footnote{Royal Society University Research Fellow and Fellow
of
Pembroke College, Cambridge}\\
{\ }\\
Department of Applied Mathematics \& Theoretical Physics\\
University of Cambridge, Cambridge CB3 9EW\\
\end{center}
\begin{center}
June 1997
\end{center}

\vspace{20pt}
\begin{quote}\baselineskip 13pt
\noindent{\bf Abstract}
The q-monopole bundle introduced previously is extended to 
a general construction for quantum group bundles with non-universal
differential calculi. We show that the theory applies to several
other classes of bundles as well, including bicrossproduct quantum
groups, the quantum double and combinatorial bundles associated to
covers of compact manifolds.
\end{quote}
\baselineskip 23pt

\section{Introduction}

A `quantum group gauge theory' in the sense of bundles with total
and base `spaces' noncommutative algebras (and quantum gauge group)
has been introduced in \cite{BrzMa:gau} with the construction of
the $q$-monopole over the $q$-sphere. Two nontrivial features of
this $q$-monopole are the  use of non-universal quantum
differential calculi and  construction in terms of patching of
trivial bundles. Several aspects of general formalism concerning
nonuniversal calculi were left open, however, and in the present
paper we study some of these aspects further, providing a
continuation of the general theory in \cite{BrzMa:gau}.

We recall that in noncommutative geometry the nonuniqueness of the
differential calculus is much more pronounced than it is
classically. Although every algebra has a universal or `free'
calculus it is much too large and one has to quotient it if one is
to have quantum geometries `deforming' the classical situation.
There are many ways to do this, however, and even for quantum
groups (where we can demand (bi)covariance) the calculus is far
from unique. In the case of quantum principal bundles one needs
quantum differential calculi both on the base and on the quantum
group fibre which have to fit together to provide a nontrivial
calculus on the total space. This is the problem which we address
here and its solution is the main result of the present paper. We
introduce in Section~3 a natural construction which builds up the
calculus on the total space of the bundle  from specified
`horizontal forms' related to the base, a specified bicovariant
calculus on the quantum group fibre and a connection form on the
bundle with the universal calculus. Roughly speaking, it is the
maximal differential calculus having the prescribed horizontal and
fibre parts and such that the connection form is differentiable.
This  approach  appears to be   different from and, we believe,
more complete than recent attempts on this problem in
\cite{Dur:dif}\cite{PflSch:dif}.

The remainder of the paper is devoted to examples and applications
of this construction. We re-examine the q-monopole in Section~4 and
verify that this example from \cite{BrzMa:gau} fits into the
general formalism.

In Section~5 we  consider a different application of the theory. We
show that the combinatorial data associated to a cover of a compact
manifold may be encoded in a discrete quantum differential calculus
over the indexing set of the cover. This demonstrates the novel
idea of doing (quantum) geometry of the combinatorics associated to
a manifold rather than the combinatorics of the classical geometry.
We show that the Czech cohomology may be recovered as the quantum
cohomology over the cover. We also consider quantum group gauge
theory over the cover as a potential source of new invariants of
manifolds. Note that classical differential calculi are not
possible over discrete sets, but nontrivial quantum ones are, i.e.
this is a natural use of quantum geometry.

In Section~6 we further apply the theory to construct
left-covariant quantum differential calculi on certain Hopf
algebras of cross product form. We regard them as trivial quantum
principal bundles and apply the results of Section~3.  Examples
include all cross product Hopf algebras such as the bicrossproduct
quantum groups in\cite{Ma:phy}, the biproducts and
bosonisations\cite{Ma:skl}\cite{Ma:poi} and the quantum double \cite{Ma:mec}.
Although the bundles here are `trivial', the uniform construction
of natural quantum differential calculi on them by abstract methods
would be a first step towards their patching to obtain nontrivial
bundles.

We begin the paper in Section~2 with some preliminaries from
\cite{BrzMa:gau}\cite{Haj:str}\cite{Brz:tra}, including the
definition of a trivial quantum principal bundle. This is basically
an algebra factorising as $P=MH$ where $M$ is the `base' algebra
and $H$ is a quantum group. All algebras in the paper should be
viewed as `coordinates' although, when the algebra is
noncommutative, they will not be the actual coordinate ring of any
usual manifold. For quantum groups, we use the notations and
conventions in \cite{Ma:book}. In particular, $\Delta:H\to H\tens
H$ denotes the coproduct expressing the `group structure' of
quantum group $H$. $S:H\to H$ denotes the antipode expressing
`group inversion', and $\eps:H\to \C$ denotes the counit,
expressing `evaluation at the group identity'. We work over $\C$.
All general constructions not involving $*$ work over a general
field just as well.

\section{Preliminaries}

In this section, we recall the basic definitions and notations to
be used throughout the paper, up to and including the definition of
a quantum principal bundle with nonuniversal calculus from
\cite{BrzMa:gau}. The same formalism has been extended to braided
group fibre and, beyond, to merely a coalgebra as fibre of the
principal bundle\cite{BrzMa:coa}, to which some of the results in
the paper should extend.

If $P$ is an algebra, we denote by $\Omega^1 P$ its universal or
K\"ahler differential structure or quantum cotangent space. Here
$\Omega^1 P=\ker\mu\subset P\tens P$, where $\mu$ is the product
map. The differential $\extd_U:P\to\Omega^1 P$ is $\extd_U u=1\tens
u-u\tens 1$. We denote by $\Omega^1(P)$ a general nonuniversal
differential structure or cotangent space. By definition, this is a
$P$-bimodule and a map $\extd:P\to \Omega^1(P)$ obeying the Leibniz
rule and such that $P\tens P\to\Omega^1(P)$ provided by $u\tens
v\mapsto u\extd v$ is surjective. It necessarily has the  form
$\Omega^1(P)=\Omega^1 P/\CN$ where $\CN\subset\Omega^1 P$ is a
subbimodule, and $\extd=\pi_{\CN}\circ\extd_U$ where $\pi_{\CN}$ is
the canonical projection. Nonuniversal calculi are in 1-1
correspondence with nonzero subbimodules $\CN$.

When $P$ is covariant under a quantum group $H$ by a (say) right
coaction $\Delta_R:P\to P\tens H$ as a comodule algebra (i.e.
$\Delta_R$ is a coaction and an algebra map), $\Omega^1(P)$ is
right covariant (in an obvious way) iff $\Delta_R(\CN)\subset
\CN\tens H$. Here $\Delta_R$ is extended as the tensor product
coaction to $P\tens P$ and restricted to $\Omega^1P$ for this
equation to make sense. We will consider only calculi on $P$ of
this form in the paper. Similar formulae hold for left covariance.

When $H$ is a Hopf algebra the coproduct $\Delta:H\to H\tens H$ can
be viewed as both a right and a left coaction of $H$ on itself by
`translation'. We will be interested throughout in nonuniversal
differential calculi $\Omega^1(H)$ which are both left and right
covariant (i.e. bicovariant) under $\Delta$. The subbimodule $\CN$
in the left covariant case in necessarily of the form
$\CN=\theta(H\tens\CQ)$ where $\theta:H\tens H\to H\tens H$ is
defined by
\eqn{theta}{\theta(g\tens h)=g Sh\o\tens h\t}
and $\CQ\subset
\ker\eps\subset H$ is a right ideal. Left covariant calculi  are in
1-1 correspondence with such $\CQ$ \cite{Wor:dif}.
Bicovariant calculi
$\Omega^1(H)$ are in 1-1 correspondence with right ideals $\CQ$
which are in addition stable under $\Ad$ in the sense
$\Ad(\CQ)\subset\CQ\tens H$ \cite{Wor:dif}. Here $\Ad$ is the right adjoint
coaction $\Ad(h)=h\t\tens (Sh\o)h\th$. We use in these formulae the
notation $\Delta h=h\o\tens h\t$ (summation understood) of the
resulting element of $H\tens H$, and higher numbers for iterated
coproducts. The universal calculus on $H$ is bicovariant and
corresponds to $\CQ=\{0\}$.

The space $\ker\eps/\CQ$ is the space of left-invariant 1-forms on
$H$. We denote by $\pi_{\CQ}$ the canonical projection. The dual of
$\ker\eps/\CQ$ (suitably defined) is the space of left-invariant
vector fields or `invariant quantum tangent space' on $H$. Hence a
map which classically has values in the Lie algebra of gauge group
will be formulated now as a map from $\ker\eps/\CQ$. This is the
approach in \cite{BrzMa:gau} for connections with nonuniversal
calculi. Note that it depends on the choice of calculus. The moduli
of bicovariant calculi (or more precisely, of quantum tangent
spaces) on a general class of quantum groups has been obtained  in
\cite{Ma:cla}; it is typically discrete but infinite.

Since a general differential calculus is the projection of a
universal one, it is natural to consider principal bundles and
gauge theory with the universal calculi $\Omega^1P$, $\Omega^1H$
first, and construct the general bundles by making quotients.
Therefore, we recall first the definitions for this universal case.
A quantum principal $H$-bundle with the universal calculus is an
$H$-covariant algebra $P$ as above, such that the map
$\chi:\Omega^1P\to P\tens \ker\eps$ defined by $\chi(u\tens
v)=u\Delta_Rv$  is surjective and obeys $\ker\chi=P(\Omega^1M)P$,
where $M=\{u\in P|\Delta_Ru=u\tens 1\}$ is the invariant
subalgebra. The latter plays the role of coordinates of the
`base'. For a complete theory, we also require that $P$ is flat as
an $M$-bimodule.  The surjectivity of $\chi$ corresponds in the
geometric case to the action being free. The kernel condition  says
that the joint kernel of all `left-invariant vector fields
generated by the action' (the maps $\Omega^1P\to P$ obtained by
evaluating against any element of $\ker\eps^*$) coincides with the
`horizontal 1-forms' $P(\Omega^1M)P$ pulled back from the base. It
plays the role in the proofs in \cite{BrzMa:gau} played classically
by local triviality and dimensional arguments. The surjectivity and
kernel conditions are equivalent to $\chi_M:P\tens_MP\to P\tens H$
being a bijection, where $\chi$ descends to the map $\chi_M$
(cf. \cite[Proposition~1.6]{Haj:str},
\cite[Lemma~3.2]{Brz:tra}). This
is the Galois condition arising independently in a more algebraic
context,  cf \cite{Mon:hop}
(not connected with
connections and differential
structures, however). We prefer to list the two conditions
separately for conceptual reasons.

A connection $\omega_U$ on a quantum principal bundle with
universal calculus is a map $\omega_U:\ker\eps\to
\Omega^1P$ such that $\chi\circ\omega_U=1\tens\id$ and
$\Delta_R\circ\omega_U=(\omega_U\tens\id)\circ\Ad$. It is shown in
\cite{BrzMa:gau} that such connections are in 1-1 correspondence
with equivariant  complements to the horizontal forms
$P(\Omega^1M)P\subset\Omega^1 P$. We are now ready for the general
case:

\begin{defin}\cite{BrzMa:gau} A general quantum principal bundle
$P(M,H,\CN,\CQ)$ is an $H$-covariant algebra $P$, an $H$-covariant
calculus $\Omega^1(P)$ described by subbimodule $\CN$ and a
bicovariant calculus $\Omega^1(H)$ described by $\Ad$-invariant
right ideal $\CQ$ compatible in the sense $\chi(\CN)\subseteq
P\tens\CQ$ and such that the map $\chi_{\CN}:\Omega^1(P)\to
P\tens\ker\eps/\CQ$ defined by $\chi_{\CN}\circ
\pi_{\CN}=(\id\tens\pi_{\CQ})\circ\chi$ is surjective and has
kernel $P(\extd M)P$.
\end{defin}

The surjectivity and kernel conditions here can also be
written as an exact sequence
\begin{equation}
0\rightarrow P(\extd M)P\rightarrow \Omega\sp 1(P)
\stackrel{\chi_{\CN}}{\rightarrow}
P\otimes\ker\epsilon/\CQ \rightarrow 0,
\end{equation}
and thus combined
into single `differential Galois' condition by noting that
$\chi_{\CN}$ descends to a map $\Omega^1(P)/P(\extd M)P\to
P\tens\ker\eps/\CQ$ and requiring this to be an isomorphism.
The condition $\chi(\CN)\subseteq P\tens
\CQ$ expresses `smoothness' of the action and is needed for
$\chi_{\CN}$ to be well-defined. In fact if $P(M,H,\CN,\CQ)$ is a
quantum principal bundle then the inclusion above implies the equality
$\chi(\CN)= P\tens
\CQ$ \cite[Corollary~1.3]{Haj:str}. On the other hand  if $P(M,H)$ is
already a
quantum principal bundle with the universal calculus then the
equality $\chi(\CN)=P\tens\CQ$ is sufficient to ensure
that $P(M,H,\CN,\CQ)$ is a quantum principal bundle with the
corresponding non-universal differential calculi. Conversely, if
$P(M,H,\CN,\CQ)$ is a quantum principal bundle with the
corresponding non-universal differential calculi then $P(M,H)$ is a
quantum principal bundle with the universal calculus if and only if
$\ker\chi\cap\CN\subseteq
P(\Omega^1 M) P\cap \CN$ \cite{Haj:pri}. Finally, a
connection on $P(M,H,\CN,\CQ)$ is a map $\omega:\ker\eps/\CQ\to
\Omega^1(P)$ such that $\chi_{\CN}\circ\omega=1\tens\id$ and
$\Delta_R\circ\omega=(\omega\tens\id)\circ\Ad$. The $\Ad$ here denotes
the quotient of the right adjoint coaction on $H$ to the space
$\ker\eps/\CQ$ given by $\Ad\circ\pi_{\CQ} =
(\pi_{\CQ}\tens\id)\circ\Ad$. As  explained in
\cite{BrzMa:gau}, connections are in 1-1 correspondence with
equivariant complements to the horizontal forms $P(\extd
M)P\subset\Omega^1(P)$. See
\cite{BrzMa:gau}\cite{Brz:tra}\cite{Ma:war95} for further details
and formalism in this approach.

There are also two main general constructions for bundles and
connections in \cite{BrzMa:gau}, the first of them used to construct
the local patches of the $q$-monopole and the second of them used
to construct the $q$-monopole globally.

\begin{example}\cite{BrzMa:gau} Let $P$ be an $H$-covariant algebra
and suppose $\Phi:H\to P$ is a convolution-invertible linear map
such that $\Phi(1)=1$ and
$\Delta_R\circ\Phi=(\Phi\tens\id)\circ\Delta$. Then $M\tens H\to P$
by $m\tens h\mapsto m\Phi(h)$ is a linear isomorphism and
$P(M,H,\Phi)$ is a quantum principal bundle with universal calculus.
There is a connection
\eqn{omegabeta}{
\omega_U(h)=\Phi^{-1}(h\o)\beta_U(\pi_\eps(h\t))\Phi(h\th)
+\Phi^{-1}(h\o)\extd_U\Phi(h\t)}for any
$\beta_U:\ker\eps\to\Omega^1M$. Here $\pi_\eps(h)=h-\eps(h)$ is the
projection to $\ker\eps$.
The case $\beta=0$ is called the {\em
trivial connection}.
\end{example}

In fact, $P$ is a cleft extension of $M$ by $H$ and has the
structure of a cocycle cross product. If, in addition, $\CQ$ and
$\CN$ define $\Omega^1(H)$ and $\Omega^1(P)$ as in Definition~2.1
then $P(M,H,\CN,\CQ)$ is a quantum principal
bundle with nonuniversal calculus. We call this a {\em  trivial
quantum principal bundle with general differential calculus}. We
will obtain in the paper the construction of connections $\omega$
from $\beta$ in this case.

\begin{example}\cite{BrzMa:gau} If $P$ is itself a Hopf algebra and
$\pi:P\to H$ a Hopf algebra surjection. $P$ becomes $H$-covariant
by $\Delta_R=(\id\tens\pi)\circ\Delta$. Suppose that the product
map $\ker\pi|_M\tens P\to \ker\pi$  is surjective. Then
$P(M,H,\pi)$ is a quantum principal bundle with universal calculus.
If there is a linear map $i:\ker\eps_H\to \ker\eps_P$ such that
$\pi\circ i=\id$ and
$(\id\tens\pi)\circ\Ad\circ i=(i\tens\id)\circ\Ad$, then there is a
connection
\eqn{omegai}{\omega_U(h)=(Si(h)\o)\extd i(h)\t}
It is
called the {\em canonical connection} associated to a linear
splitting $i$.
\end{example}
\begin{remark}
\rm Note that if $P$ and $H$ are Hopf algebras and
$\pi:P\to H$ is a Hopf algebra surjection, then the canonical map
$\chi$ is surjective since   it is obtained by projecting the
inverse $\theta^{-1}$ of the linear automorphism $\theta$ of  $P\tens
P$  in (\ref{theta}) down to $P\tens
H$, i.e. $\chi = (\id\tens\pi)\circ\theta^{-1}$. The condition that
the product map $\ker\pi|_M\tens P\to \ker\pi$  be  surjective
provides that the kernel of $\chi$ is equal to horizontal
one-forms. Combining \cite[Theorem~I]{Sch:pri} with
\cite[Lemma~1.3]{Sch:nor} one finds that $\ker\pi|_M\tens P\to
\ker\pi$  is  surjective if there is a linear map $j :H\to P$ such
that $j(1)=1$ and $\Delta_R\circ j = (j\tens\id)\circ\Delta$. More
precisely, \cite[Theorem~I]{Sch:pri} and
\cite[Lemma~1.3]{Sch:nor} imply that if such a $j$ exists then in
addition to $P(M,H,\pi)$ there is also a quantum principal bundle
$P(M,H',\pi')$ where $H' = P/(\ker\pi|_M\cdot P)$. Therefore one can write
the following commutative diagram
$$
\begin{CD}
@. 0 @>>>    0 @>>>  \ker s @.\\
  @.      @VVV   @VVV @VVV\\
0 @>>>    P(\Omega^1M)P   @>>>    P\tens P
  @>(\id\tens\pi')\circ\theta^{-1}>> P\tens H'@>>> 0 \\
  @.      @VVV   @VVV @VVsV\\
0 @>>>    P(\Omega^1M)P   @>>>    P\tens P
  @>(\id\tens\pi)\circ\theta^{-1}>> P\tens H @>>> 0 \\
@.     @VVV   @VVV @VVV\\
@. 0 @>>>    0    @>>> {\rm coker}s @.
\end{CD}
$$
The second and third row are exact by definition of a quantum
principal bundle. Obviously ${\rm coker}s =0$. The application of the
snake lemma (cf. \cite[Section~1.2]{Bou:hom})  yields $\ker
s =0$, i.e. $H'\subseteq H$. Since $H=P/\ker\pi$, $H'=P/(\ker\pi|_MP)$
this implies that the product map $\ker\pi|_M\tens P\to \ker\pi$
is  surjective as required.

If there exists a left integral on $H$ , i.e. $\lambda \in H^*$
such that $\lambda(1)$ = 1 and $(\lambda\tens\id)\circ\Delta =
\lambda$ then the map $j:H\to P$ can be defined by
$j = \eta_P\circ\lambda$, where $\eta_P:\C\to P$ is the unit map. The
map $j$ is clearly an intertwiner since
$$
\Delta_Rj(h) = \lambda(h) 1\tens 1 = \lambda(h\o)\tens h\t =
(j\tens\id)\Delta(h).
$$
 In particular, if $H$ is a compact quantum
group in the sense of \cite{Wor:com} then $\lambda$ is
the Haar measure on $H$.
Therefore if $H$ is a compact quantum group then the Hopf algebra
surjection $\pi: P\to H$ leads immediately to the bundle
$P(M,H,\pi)$. This fact is also proven directly by using
representation theory of compact quantum groups in \cite{Dur:hop}.
\label{remark.intertwiner}
\end{remark}

In the situation of Example~2.3, if  $\Omega^1(P)$ is left covariant with its
corresponding right ideal $\CQ_P\subseteq\ker\eps\subset P$ obeying
$(\id\tens\pi)\circ\Ad(\CQ_P)\subset \CQ_P\tens H$, then
$\Omega^1(H)$ defined by $\CQ=\pi(\CQ_P)$ provides a quantum
principal bundle $P(M,H,\pi,\CQ_P)$. We call it a {\em homogeneous
space bundle with general differential calculus}.
If $i:\ker\eps_H\to \ker\eps_P$ is
as above and, in addition, $i(\CQ)\subset\CQ_P$ then
$\omega(h)=(Si(h)\o)\extd i(h)\t$ is a connection. A refinement of
this construction
 will be provided in the paper.

\section{Differential Calculi on Quantum Principal Bundles}

In this section we obtain the main tool in the paper. This is a new
construction for general quantum principal bundles with
nonuniversal calculi, starting with a specified bicovariant calculus
$\Omega^1(H)$ on the fibre and a specified `horizontal calculus' on
the base. In the classical case one has local triviality and one
accordingly takes the calculus on $P$ coinciding with its direct
product form over each open set. That  this is actually the
standard calculus on $P$ is consequence of the smoothness part of
the axiom of local triviality. This is our motivation now.

As recalled in the Preliminaries, in the quantum case we actually
have global conditions playing the role of local
triviality\cite{BrzMa:gau}, which is the `global approach' which we
describe first. Building up the calculus on $P$ globally in this
way means  that we construct $\Omega^1(P)$ as the direct sum of a
part from the base and a part from the fibre, i.e. actually the
same process as building a connection $\omega$. Therefore, the
nonuniversal bundles constructed in this way will automatically
have the property of existence of a natural connection.

On the other hand, the data going into the construction of
$\Omega^1(P)$ should not already assume the existence of a bundle,
as this is to be constructed. Instead, the additional input data
besides the desired calculi $\Omega^1(H)$ and on the base should be
related to the `topological' and not `differential' splitting. We
therefore take for this additional `gluing' datum a connection
$\omega_U$ on $P$ as a quantum principal bundle with the universal
calculus.

Accordingly, we let $P(M,H)$ be a quantum principal bundle with the
universal calculus and $\Omega^1(H)$ a choice of bicovariant
calculus on $H$ defined by $\CQ\subseteq\ker\eps\subset H$. As far
as the differential calculus on the base is concerned, we can
specify $\Omega^1(M)$ by $\CN_M\subset\Omega^1M$ as an
$M$-subbimodule. More natural (and slightly more general) is to
specify a `horizontal' subbimodule $\CN_{\rm hor}$.

\begin{lemma}  Let $P(M,H)$ be a quantum principal bundle with the
universal calculus and let $\omega_U$ be a connection on it. Let $\CQ$ specify
a bicovariant calculus on $\Omega^1(H)$. Let $h_\omega : P\tens Q\tens
P \to \Omega^1P$ be a linear map given by $h_\omega(u,q,v) =
uv\bo\omega_U(qv\bt)-u\omega_U(q)v$, where we write $\Delta_Ru=u\bo\tens
u\bt$ (summation understood). Then $\CN_0={\rm Im} h_\omega
\subset
P(\Omega^1M)P$
is a $P$-subbimodule invariant under $\Delta_R$ in the sense
$ \Delta_R\CN_0\subset\CN_0\tens H.$
\label{lemma.n0}
\end{lemma}
\proof Clearly $wh_\omega(u,q,v) = h_\omega(wu,q,v)$, for any $u,v,w\in P$
and $q\in \CQ$. Also
\begin{eqnarray*}
h_\omega(u,q,v)w & = &uv\bo\omega_U(qv\bt)w-u\omega_U(q)vw\\
& = &
uv\bo\omega_U(qv\bt)w - uv\bo w\bo\omega_U(qv\bt w\bt)
+h_\omega(u,q,vw)\\
& = & h_\omega(uv\bo,qv\bt,w) + h_\omega(u,q,vw).
\end{eqnarray*}
Therefore $\CN_0 = {\rm Im} h_\omega$ is a subbimodule of
$\Omega^1P$. Furthermore
$$
\chi(h_\omega(u,q,v)) = \chi(uv\bo\omega_U(qv\bt)-u\omega_U(q)v) =
uv\bo\tens qv\bt - (u\tens q)(v\bo\tens v\bt) = 0,
$$
i.e. $\CN_0\in\ker\chi = P(\Omega^1M)P$. Finally,
\begin{eqnarray*}
\Delta_R(h_\omega(u,q,v))
 & = & u\bo v\bo \omega_U(q\t
v\bt\th)\tens u\bt v\bt\o Sv\bt\t Sq\o q\th v\bt\fo\\
&& - u\bo v\bo \omega_U(q\t)\tens u\bt  Sq\o q\th v\bt \\
& = & u\bo v\bo \omega_U(q\t
v\bt\o)\tens u\bt Sq\o q\th v\bt\t\\
&&  - u\bo v\bo \omega_U(q\t)\tens u\bt  Sq\o q\th v\bt\\
& = & h_\omega(u\bo, q\t, v\bo)\tens u\bt Sq\o q\th v\bt \in
\CN_0\tens H
\end{eqnarray*}
where we used the covariance of $\omega_U$ and the fact that $\CQ$ is
Ad-invariant. Therefore $\CN_0$ is right-invariant as stated.
\endproof

We call any $\Delta_R$-invariant $P$-subbimodule of
$P(\Omega^1M)P$ `horizontal'. We fix one, denoted $\CN_{\rm hor}$.
The corresponding quotient $\Omega^1_{\rm
hor}=P(\Omega^1M)P/\CN_{\rm hor}$ is our choice of `horizontal'
part of the desired calculus on $P$.

\begin{theorem} Let $P(M,H)$, $\omega_U$ be a quantum principal bundle with
the
universal calculus and connection as above. Let
$\CQ$ specify $\Omega^1(H)$ and
\[ \CN_0\subseteq \CN_{\rm hor} \subseteq P(\Omega^1M)P\]
specify $\Omega^1_{\rm
hor}$. Then
\[ \CN=\<\CN_{\rm hor},P\omega_U(\CQ)P\>\]
specifies a differential calculus $\Omega^1(P)$ with the property that
$P(M,H,\CN,\CQ)$ is
a quantum principal bundle, $P(\extd M)P=\Omega^1_{\rm hor}$ and
$\omega=\pi_{\CN}\circ\omega_U$ is a connection on the bundle.

The calculus resulting from the choice $\CN_{\rm hor}=\CN_0$ is called
the {\em maximal}
differential calculus compatible with $\omega_U$. The choice $\CN_{\rm
hor}=P(\Omega^1M)P$ is called the {\em minimal} differential calculus
compatible 
with $\omega_U$.
\label{theorem.calculus}
\end{theorem}
\proof By assumption, $\Delta_R\CN_{\rm
hor}\subset\CN_{\rm hor}\tens H$. Also
\align{\Delta_R(u\omega_U(q)v)\equad&&=u\bo \omega_U(q)\bo v\bo
\tens u\bt\omega_U(q)\bt v\bt=u\bo\omega_U(q\o)v\bo\tens
u\bt(Sq\o)q\th v\bt}
for all $u,v\in P$ and $q\in H$.  The result is manifestly in
$P\omega_U(\CQ)P\tens H$ since $\CQ$ is $\Ad$-invariant. Hence
$\Delta_R\CN\subset\CN\tens H$.

Next, we clearly have $\chi(\CN_{\rm hor})=0$. Then
$\chi(u\omega_U(q)v)=u\chi(\omega_U(q))\Delta_Rv
=uv\bo\tens qv\bt\in P\tens \CQ$
since $\CQ$ is a right ideal. Conversely, if $u\tens q\in P\tens
\CQ$ then $u\omega_U(q)\in P\omega_U(Q)P$ and
$\chi(u\omega_U(q))=u\tens q$. Hence $\chi(\CN)=P\tens Q$. We
therefore have quantum principal bundle with nonuniversal
differential calculus $\Omega^1(P)$.

Clearly $\pi_{\CN}\circ\omega_U(\CQ)=0$, hence this descends to a
map $\omega:\ker\eps/\CQ\to \Omega^1(P)$. Moreover,
\align{\chi_{\CN}\circ\pi_{\CN}\circ\omega_U\equad&&
=(\id\tens\pi_{\CQ})\circ\chi\circ\omega_U=1\tens\pi_{\CQ}}
on $\ker\eps$, where the first equality is the definition of
$\chi_{\CN}$ and the second is the equivariance of $\omega_U$.
Also,
\[\Delta_R\circ\pi_{\CN}\circ\omega_U=(\pi_{\CN}\tens\id)
\circ\Delta_R\omega_U=(\pi_{\CN}\circ\omega_U\tens\id)\circ\Ad.\]
The first equality is clear from the definition of $\Delta_R$.
Hence we have a connection $\omega$.

Finally, we note that the stated $\Omega^1(P)$ is uniquely determined by
$\omega_U$
and $\CN_{\rm hor}$ as the universal calculus with the stated properties.
Thus, suppose that $\CN'$ defines another quantum
differential calculus on $P$ such that $\pi_{\CN'}\circ\omega_U$ is
a connection. Then $\omega_U(\CQ)\subset\CN'$. The stated $\CN$ is
clearly the minimal subbimodule containing $P\omega_U(\CQ)P$ and
$\CN_{\rm hor}$, i.e. any other such $\Omega^1(P)$ is a quotient. \endproof

There is a natural generalisation of this theorem in which we
assume only that $P$ is an $H$-comodule algebra (i.e. without going
through the assumption that $P(M,H)$ is already a quantum principal
bundle with the universal calculus). For this version we assume the
existence of an $\Ad$-equivariant map
$\omega_U:\ker\eps\to\Omega^1P$ obeying
$\chi\circ\omega_U=1\tens\id$ and $\CQ$, $\CN_{\rm hor}$ as above.
Then the  map $\chi_{\CN}$ can be defined and
if its kernel is $P(\extd M) P$ then the same
conclusion holds.

We now consider how our construction looks for the two examples of
quantum principal bundles with the universal calculus in the
Preliminaries section.

\begin{propos} Consider a trivial quantum principal bundle
$P(M,H,\Phi)$
with the universal calculus and let $\Omega^1(H)$ and $\Omega^1(M)$
be determined by $\CQ$ and $\CN_M$. Then for any
$\beta_U:\ker\eps\to\Omega^1 M$
there is a differential calculus $\Omega^1(P)$ with $\Omega^1_{\rm hor}
=P(\extd M)P$ and  forming a trivial quantum principal bundle, and
\begin{equation}
\omega(h)=\Phi^{-1}(h\o)\beta\circ\pi_\eps(h\t)\Phi(h\th)
+\Phi^{-1}(h\o)\extd\Phi(h\t)
\label{strong.beta}
\end{equation}
is a connection on it for $\beta:\ker\eps \to \Omega^1_P(M)$, where  obtained
by restricting $\Omega^1_P(M)= \pi_{\CN}(\Omega^1M)$.
\label{prop.triv}
\end{propos}
\proof  We define $\omega_U :\ker\eps \to \Omega^1P$ by 
$\omega_U(h)=\Phi^{-1}(h\o)\beta_U(\pi_\eps(h\t))\Phi(h\th)
+\Phi^{-1}(h\o)\extd_U\Phi(h\t)$ as a connection on the bundle with
universal calculus.  We also take $\CN_{\rm hor}=\<P\CN_MP,\CN_0\>$ where
$\CN_0$ is determined
by $\omega_U$. We can now apply Theorem~3.2. Note also that $\CN =
\<P\CN_MP, P\omega_U(\CQ)P\>$ as $\CN_0\subset P\omega_U(\CQ)P$.

Explicitly, the sub-bimodule corresponding to $\Omega^1(P)$ is
$\CN=\<\CN_{\rm hor}, P\hat\Phi(\CN_H)P\>$ where
$\CN_H=\theta(H\tens \CQ)$ and
$\hat\Phi:H\tens H\to\Omega^1P$, 
\[\hat\Phi(g\tens h)
=\Phi(gh\o)\Phi^{-1}(h\t)\beta_U(h\th)\Phi(h\fo)+\Phi(gh\o)
\Phi^{-1}(h\t)
\tens\Phi(h\th)\]
for all $h,g\in H$. This makes it clear that we recover here the
construction for nonuniversal trivial bundles in terms of a map
$\hat\Phi$ in \cite{Brz:the}. Note that $\hat\Phi\circ\theta(g\tens
h)=\hat\Phi(gSh\o\tens h\t)
=\Phi(g)\omega_U(h)$, for any $g,h\in H$.  Note that the inherited
differential structure on $M$, $\Omega^1_P(M)\subset \Omega^1(P)$ is
smaller than the original $\Omega^1(M) =\Omega^1M/\CN_M$ unless 
$\CN_0\cap \Omega^1M\subseteq\CN_P$. \endproof

{}From the proof of Theorem~3.2 we see that the resulting trivial
bundle with nonuniversal calculus is of the general type discussed
after Example~2.2; we succeed by the above to put a general class
of connections on it. Also note that we may take more general
$\CN_{\rm hor}$ and any $\beta_U:\ker\eps\to\Omega^1M$ to arrive at
some $\Omega^1(P),\omega$, though not necessarily of the form
stated.

\begin{lemma}
Let $P(M,H,\Phi,\CQ,\CN)$ be a trivial quantum principal bundle with
differential calculus determined by $\CQ$ and $\CN$. Then
$\beta:\ker\eps\to\Omega^1_P(M)$ defines a connection $\omega$ by
(\ref{strong.beta}) in
Proposition~3.3. if and only if for all $q\in \CQ$,
\begin{equation}
\Phi^{-1}(q\o)\beta(\pi_\eps(q\t))\Phi(q\th) =
-\Phi^{-1}(q\o)\extd\Phi (q\t). 
\label{cond.beta}
\end{equation}
Furthermore, if $\Phi$ is an algebra map then for all $h\in H$,
\[\Phi^{-1}(q\o)\beta(\pi_\eps(q\t h))\Phi(q\th) =
\eps(h)\Phi^{-1}(q\o)\beta(\pi_\eps(q\t))\Phi(q\th) .\]
\label{lemma.cond.beta}
\end{lemma}
\proof Requirement (\ref{cond.beta}) is another way of expressing the
fact that $\omega(q) = 0$ for all $q\in \CQ$. Since $\CQ$ is a right
ideal, condition (\ref{cond.beta}) implies that
\[ h\o\tens\Phi^{-1}(q\o h\t)\beta(\pi_\eps(q\t h\th))\Phi(q\th
h\fo)\tens h\fiv =
-h\o\tens \Phi^{-1}(q\o h\t)\extd \Phi(q\t h\th)\tens h\fo.
\]
Applying $\Phi\tens\id\tens\Phi^{-1}$ and multiplying we thus obtain
\begin{eqnarray*}
\Phi(h\o)\Phi^{-1}(q\o h\t)\beta(\pi_\eps(q\t h\th))\Phi(q\th
h\fo)\Phi^{-1}(h\fiv) & = &\\
&&\!\!\!\!\!\!\!\!\!\!\!\!\!\!\!\!\!\!-\Phi(h\o)\Phi^{-1}(q\o
h\t)\extd \Phi(q\t h\th)\Phi^{-1}(h\fo).
\end{eqnarray*}
If $\Phi$ is an algebra map the above formula simplifies further
\[\Phi^{-1}(q\o)\beta(\pi_\eps(q\t h))\Phi(q\th)=
-\Phi^{-1}(q\o)\extd (\Phi(q\t)\Phi( h\o))\Phi^{-1}(h\t).
\]
The application of the Leibniz rule and the fact that $q\in\ker\eps$
the yields
\[\Phi^{-1}(q\o)\beta(\pi_\eps(q\t h))\Phi(q\th
h)=
-\eps(h)\Phi^{-1}(q\o)\extd \Phi(q\t),
\]
which in view of (\ref{cond.beta}) implies the assertion. Notice that
the condition one obtains in this way deals entirely with the structure
of $\Omega_{\rm hor}$ and is the consequence of the existence of
$\CN_0$. \endproof

\begin{propos} Consider a quantum principal bundle with the universal
calculus of the homogeneous type $P(M,H,\pi)$ where $\pi:P\to H$ is
a Hopf algebra surjection. For any $\Omega^1(H)$, if $\omega_U$
is left-invariant and $\CN_{\rm hor}$ is left-invariant under the
left-regular coaction of $P$ as a Hopf algebra, then
$\Omega^1(P)$ in Theorem~3.2 is left covariant. Moreover,
left-invariant $\omega_U$
are canonical connections in 1-1 correspondence with $i$ as in
Example~2.3. Left-covariant $\CN_{\rm hor}$ are in 1-1 correspondence with
right
ideals $\CQ_0\subset\CQ_{\rm hor}\subseteq\ker\pi$, where
\[ \CQ_0=\span \{i(q)u-i(q\pi(u))| \ q\in\CQ,\ u\in P\}.\]
\label{prop.homog}
\end{propos}
\proof We can regard $\CN_{\rm hor}$ as a
subbimodule of $\Omega^1P$.
As such, it defines a differential calculus $\Omega^1P/\CN_{\rm
hor}$ on $P$. As $P$ is now a Hopf algebra, the calculus is left
covariant i.e. $\Delta_L\CN_{\rm hor}\subset \CN_{\rm hor}\tens P$
iff $\CN_{\rm hor}=\theta(P\tens\CQ_{\rm hor})$ for a right ideal
$\CQ_{\rm hor}\subseteq\ker\eps\subset P$. Here $\Delta_L$ is the
left regular coaction or `translation' on $P\tens P$ obtained from
the coproduct.

On the other hand, since $\CN_{\rm hor}\subset P(\Omega^1M)P$, we
know that $\chi(\CN_{\rm hor})=0$. Take $q\in\CQ_{\rm hor}$. Then
$0=\chi\theta(1\tens q)=(Sq\o)q\t\tens\pi(q\th)=\pi(q)$ so
$\CQ_{\rm hor}\subseteq\ker\pi$. Conversely, if  $\CQ_{\rm
hor}\subseteq\ker\pi$ then clearly $\CN_{\rm
hor}=\theta(P\tens\CQ_{\rm hor})\subset P(\Omega^1M)P$, since $\ker\pi
= (\ker\eps\mid_M)P$.

If the connection $\omega_U$ is invariant in the sense $\Delta_L\omega_U(h)
=1\tens\omega_U(h)$ for any $h\in \ker\eps$ then clearly $\Delta_L(u\omega_U(q)v)\in
P\tens P\omega_U(\CQ)P$ for all $u,v\in P$ and
$q\in\CQ\subset\ker\eps\subset H$. Therefore $\CN$ defined in
Theorem~3.2 obeys $\Delta_L\CN\subset P\tens\CN$, i.e.
$\Omega^1(P)$ is left covariant.

The canonical connection associated to $i$ as in Example~2.3 is
invariant:
\align{\Delta_L\omega_U(h)\equad &&= \Delta_L((Si(h)\o)\extd_U
i(h)\t)=(Si(h)\t)i(h)\th\tens Si(h)\o\tens i(h)\fo\\ &&=1\tens Si(h)\o\tens
i(h)\t=1\tens\omega_U(h).} Conversely if
$\omega_U$ is an invariant connection , we define
$i:\ker\eps_H\to\ker\eps_P$ by $i =
(\eps_P\tens\id)\circ\omega_U$. Let $\theta^{-1}$ be the inverse to
the canonical map 
$\theta : P\tens P\to P\tens P$ defined as in (\ref{theta}). Explicitly
$\theta^{-1}(u\tens v) =
uv\o\tens v\t$. Clearly $\theta^{-1} =
(\id\tens\eps_P\tens \id)\circ\Delta_L$. Since $\omega$ is left-invariant
one immediately finds that $\theta^{-1}\circ\omega_U(h) = 1\tens
i(h)$. Thus $\omega_U(h) = \theta(1\tens i(h)) = Si(h)\o\tens
i(h)\t = Si(h)\o\extd_Ui(h)\t$ and $\omega_U$ has the structure of the
canonical connection associated to $i$. It remains to prove that $i$
is an $\Ad$-covariant splitting. Since $\chi =
(\id\tens\pi)\circ\theta^{-1}$ the fact that
$\chi\circ\omega_U(h) = 1\tens h$ implies that $\pi(i(h)) = h$. Finally
compute
$$
\Delta_R(\omega_U(h)) = Si(h)\t\tens i(h)\th \tens \pi(Si(h)\o i(h)\fo).
$$
On the other hand $\omega_U$ is a connection therefore
$$
\Delta_R(\omega_U(h)) = \omega_U(h\t)\tens Sh\o h\th = Si(h\t)\o \tens
i(h\t)\t \tens  Sh\o h\th.
$$
Applying $(\eps_P\tens\id\tens\id)$ to above equality one obtains the
required $\Ad$-covariance of $i$.

Using the fact that $\omega_U$ is left-invariant we find
\align{
\Delta_L(h_\omega(u,q,v)) &= &u\o v\o\tens u\t v\t\omega_U(q\pi(v\th)) -
u\o v\o\tens u\t v\t\omega_U(q)\\& = &u\o v\o \tens h_\omega(u\t,q,v\t),}
where $h_\omega$ is the map defined in Lemma~\ref{lemma.n0}. Therefore
$\CN_0$ is left-invariant and there is corresponding right ideal
$\CQ_0 \in \ker\eps_P$ given by $\CN_0= \theta(P\tens \CQ_0)$. Since
$\CN_{\rm hor}$ contains
necessarily $\CN_0$, the right ideal
$\CQ_{\rm hor}$ must contain $\CQ_0$. For the
canonical connection induced by the splitting $i$, $\CQ_0$ comes out as
stated. The fact that $\CQ_0$ is a right ideal
can be established directly since
$$
(i(q)u-i(q\pi(u)))v = (i(q)uv-i(q\pi(uv))) -
(i(q\pi(u))v-i(q\pi(u)\pi(v))) \in \CQ_0.
$$

\begin{sloppypar}
For completeness, we also show that the resulting bundle  is indeed
of the natural nonuniversal homogeneous type discussed after
Remark~\ref{remark.intertwiner}. First of all   note  that
$\theta^{-1}(u\omega_U(q)v)=u(Si(q)\o)i(q)\t
v\o\tens i(q)\th v\t=uv\o\tens i(q)v\t$. Hence
$\CN=\theta(P\tens\CQ_P)$ where $\CQ_P=\<\CQ_{\rm hor},i(Q)P\>$.
{}From this it is also clear that $\Omega^1(P)$ is left covariant,
as $\CQ_P$ is clearly a right ideal. Also,  $\pi(\CQ_P) =
\CQ$. It remains to verify whether $(\id\tens\pi)\Ad (\CQ_P)
\subset \CQ_P\tens H$. Take any $q\in
\CQ_P$, then $Sq\s 1
\tens q\s 2 \in \CN$. By construction $\Omega^1(P)$ is right
$H$-covariant, therefore $Sq\s 2\tens q\s 3\tens \pi(Sq\s 1 q\s 4)
\in \CN\tens H$. Applying $\theta^{-1}\tens \id$ to this one thus
obtains that $1\otimes q\s 2 \tens \pi(Sq\s 1 q\s 3)\in P\otimes
\CQ_P\otimes H$. Therefore $(\id\tens\pi)\Ad
(\CQ_P) \subset \CQ_P\tens H$ as required.
\end{sloppypar}

 In this case it is clear that
$i(\CQ)\subset \CQ_P$, i.e. the canonical connection is of the type
mentioned after Remark~2.4 from \cite{BrzMa:gau}.
\endproof

In the case of a homogeneous quantum principal bundle with a general
differential calculus of the type  discussed after
Remark~2.4. we can establish the one-to-one  correspondence between
invariant connections and $\Ad$-covariant splittings as follows.
The conditions
satisfied by $\CQ_P$ and $\CQ$ allow for definition of maps
$\overline{\Ad} : \ker\eps_P/\CQ_P \to \ker\eps_P/\CQ_P\tens H$ and
$\overline{\pi} : \ker\eps_P/\CQ_P \to \ker\eps_H/\CQ$ by
$\overline{\Ad}\circ\pi_{\CQ_P} = (\pi_{\CQ_P}\tens\pi)\circ\Ad$ and
$\overline{\pi}\circ\pi_{\CQ_P} = \pi_{\CQ}\circ\pi$. Here $\pi_{\CQ}
:\ker\eps_H\to
\ker\eps_H/\CQ$ and $\pi_{\CQ_P} :\ker\eps_P\to
\ker\eps_P/\CQ_P$  are canonical surjections.

\begin{propos}
The left-covariant connections $\omega$ in $P(M,H,\pi,\CQ_P)$ are in
one-to-one correspondence with the linear maps
$i:\ker\eps_H/\CQ\to\ker\eps_P/\CQ_P$ such that
$\overline{\pi}\circ i=\id$ and
$\overline{\Ad}\circ i=(i\tens\id)\circ\Ad$.
\end{propos}
\proof Assume that $\omega:\ker\eps_H/\CQ\to\Omega^1(P)$ is an
invariant connection in $P(M,H,\pi,\CQ_P)$. Define a map
$\overline{\eps}:\Omega^1(P)\to\ker\eps_P/\CQ_P$ by the commutative
diagram with exact rows
$$
\begin{CD}
0 @>>>    \CN    @>>>    \Omega^1P @>{\pi_{\CN}}>> \Omega^1(P) @>>> 0\\
  @.      @VV{\eps_P\tens\id}V   @VV{\eps_P\tens\id}V @VV{\overline{\eps}}V\\
0 @>>>    \CQ_P    @>>>    \ker\eps_P @>{\pi_{\CQ_P}}>> \ker\eps_P/\CQ_P @>>>
  0
\end{CD}
$$
Let $i = \overline{\eps}\circ\omega$. Then we have the following
commutative diagram with exact rows
\begin{equation}
\begin{CD}
0 @>>>    \CN    @>>>    \Omega^1P @>{\pi_{\CN}}>> \Omega^1(P) @.\\
  @.      @VV{\Delta_L}V   @VV{\Delta_L}V @VV{\Delta_L}V\\
0 @>>> P\tens \CN  @>>> P\tens\Omega^1P @>{\id\tens \pi_{\CN}}>>
  P\tens\Omega^1(P) @>>> 0\\
  @.      @VV{\id\tens\eps_P\tens\id}V   @VV{\id\tens\eps_P\tens\id}V
  @VV{\id\tens\overline{\eps}}V\\
0 @>>> P\tens\CQ_P    @>>> P\tens\ker\eps_P @>{\id\tens\pi_{\CQ_P}}>> P\tens
  \ker\eps_P/\CQ_P @>>> 0\\
@.      @VV{\id\tens\pi}V   @VV{\id\tens\pi}V
  @VV{\id\tens\overline{\pi}}V\\
@. P\tens\CQ    @>>> P\tens\ker\eps_H @>{\id\tens\pi_{\CQ}}>> P\tens
  \ker\eps_H/\CQ @>>> 0
\end{CD}
\label{diag.chi}
\end{equation}
The first two maps $\Delta_L$ are left coactions of $P$ on $P\tens P$
obtained from the left regular coaction of $P$ provided by the
coproduct  while the third $\Delta_L$ is their
projection to $\Omega^1(P)$. The third column gives the map $\chi$,
therefore the fourth column describes $\chi_{\CN}$, i.e. $\chi_{\CN}
=
(\id\tens\overline{\pi}\circ\overline{\eps})\circ\Delta_L$.
Since $\omega$ is a connection, $\chi_{\CN}(\omega(h)) = 1\tens h$ for
any $h\in \ker\eps_H/\CQ$. Thus we have
$$
1\tens h =(\id\tens\overline{\pi}\circ\overline{\eps})\circ\Delta_L(\omega(h))
= 1\tens\overline{\pi}\circ\overline{\eps}\circ\omega(h) = 1\tens
\overline{\pi}\circ i(h).
$$
To derive the second equality we used invariance of
$\omega$. Therefore $\overline{\pi}\circ i =\id$.

Next consider the map $\theta_N :
\Omega^1(P)\to P\tens\ker\eps_P/\CQ_P$, given by $\theta_N =
(\id\tens\overline{\eps})\circ\Delta_L$. This map makes the following
diagram commute
$$
\begin{CD}
@. 0 @>>>    0 @>>>  \ker\theta_N @.\\
  @.      @VVV   @VVV @VVV\\
0 @>>>    \CN    @>>>    \Omega^1P @>{\pi_{\CN}}>> \Omega^1(P)@>>> 0 \\
  @.      @VV{\theta^{-1}}V   @VV{\theta^{-1}}V @VV{\theta_N}V\\
0 @>>> P\tens\CQ_P    @>>> P\tens\ker\eps_P @>{\id\tens\pi_{\CQ_P}}>> P\tens
  \ker\eps_P/\CQ_P @>>> 0\\
@.     @VVV   @VVV @VVV\\
@. 0 @>>>    0    @>>> {\rm coker}\theta_N @.
\end{CD}
$$
This diagram  is a combination of the first three rows of
(\ref{diag.chi}). Clearly ${\rm coker}\theta_N = 0$. By the snake
lemma (cf. \cite[Section~1.2]{Bou:hom}),  $\ker\theta_N =0$. Therefore $\theta_N$ is a bijection. The
left-invariance of $\omega$ implies that
$$
\theta_N(\omega(h)) =
(\id\tens\overline{\eps})\circ\Delta_L(\omega(h)) =
1\tens\overline{\eps}(\omega(h)) = 1\tens i(h)
$$
for any $h\in \ker\eps_H/\CQ$. Therefore $\omega(h) =
\theta^{-1}_N(1\tens i(h))$.

Using  $\overline{\Ad}$ and $\Delta_R$ one constructs the tensor
product coaction $\overline{\Delta}_R :
P\tens\ker\eps_P/\CQ_P\to P\tens\ker\eps_P/\CQ_P\tens H$. Then $\theta_N$
is a right $H$-comodule isomorphism. This follows from the fact that
$\theta^{-1}$ is a corresponding $H$-comodule isomorphism. Explicitly
\begin{eqnarray*}
\Delta_R(\theta^{-1}(u\tens v)) & = & \Delta_R(uv\o\tens v\t) = u\o v\o\tens
v\fo\tens\pi(u\t v\t Sv\th v\fiv)\\
& = &u\o v\o\tens
v\t\tens\pi(u\t v\th ).
\end{eqnarray*}
$\Delta_R$ here is a right coaction of $H$ on $P\tens\ker\eps_P$ built
with $(\id\tens\pi)\circ\Delta$ on $P$ and $(\id\tens\pi)\circ\Ad$ on
$\ker\eps_P$. On the other hand
$$
(\theta^{-1}\tens\id)(\Delta_R(u\tens v)) = \theta^{-1}(u\o\tens v\o)\tens
\pi(u\t v\t) = u\o v\o\tens
v\t\tens\pi(u\t v\th ),
$$
where $\Delta_R$ is a standard tensor product coaction of $H$ on
$\Omega^1P\subset P\tens P$. Therefore
\begin{eqnarray*}
\overline{\Delta}_R\circ\theta_N\circ\pi_{\CN} & = &
\overline{\Delta}_R\circ(\id\tens\pi_{\CQ_P})\circ\theta^{-1} =
(\id\tens\pi_{\CQ_P}\tens\id)\circ\Delta_R\circ\theta^{-1} \\
& = & (\id\tens\pi_{\CQ_P}\tens\id)\circ(\theta^{-1}\tens\id)\circ\Delta_R =
(\theta_N\tens\id)\circ(\pi_{\CN}\tens\id)\circ\Delta_R \\
& = & (\theta_N\tens\id)\circ\Delta_R\circ\pi_{\CN}.
\end{eqnarray*}
Therefore $\theta_N$ is an intertwiner as stated. Its inverse is also
an intertwiner. We compute
$$
\Delta_R\circ\theta^{-1}_N(1\tens i(h)) =
(\theta^{-1}_N\tens\id)\circ\overline{\Delta}_R (1\tens i(h)) =
(\theta^{-1}_N\tens\id)(1\tens\overline{\Ad}(i(h))).
$$
On the other hand, since $\omega$ is a connection  this is equal to
$\Delta_R(\omega(h)) = (\omega\otimes \id)\circ \Ad (h)$. Applying
$(\theta_N\tens\id)$ to both sides one obtains
$$
1\tens((i\tens\id)\circ\Ad(h)) =  1\tens \overline{\Ad}\circ i(h),
$$
i.e.
$\overline{\Ad}\circ i=(i\tens\id)\circ\Ad$, as required.

Conversely, given $i :\ker\eps_H/\CQ \to \ker\eps_P/\CQ_P$ with the
properties described in the proposition, one defines a map $\omega
:\ker\eps_H/\CQ \to \Omega^1(P)$ by $\omega(h) = \theta^{-1}_N(1\tens
i(h))$. The $\Ad$-covariance of $i$ implies the $\Ad$-covariance of
$\omega$ since $\theta_N^{-1}$ is an intertwiner of $\Delta_R$ and
$\overline{\Delta}_R$.  Furthermore, since $\chi_{\CN} =
(\id\tens\overline{\pi})\circ\theta_N$ from diagram (\ref{diag.chi}),
$$
\chi_{\CN}(\omega(h)) = 1\tens\overline\pi(i(h)) = 1\tens h.
$$
Therefore $\omega$ is a connection. The fact that $\omega$ obtained in
this way is left-covariant is
well-known from the theory of
left-covariant calculi \cite{Wor:dif} but we include the proof for the
completeness. First consider any $u\tens v \in P\tens P$ and
compute
$$
\Delta_L(\theta(u\tens v)) = \Delta_L(uSv\o\tens v\t) = u\o Sv\t
v\th\tens u\t Sv\o \tens v\fo = u\o\tens\theta(u\t\tens v).
$$
This implies that
$$
\Delta_L(\theta^{-1}_N(u\tens v)) = u\o\tens\theta^{-1}_N(u\t\tens v).
$$
for any $u\in P$ and $v\in\ker\eps_P/\CQ_P$. Therefore
$$
\Delta_L\omega(h) = \Delta_L\circ\theta^{-1}_N(1\tens i(h)) = 1\tens
\theta_N^{-1}(1\tens i(h)) = 1\tens \omega(h),
$$
for any $h\in \ker\eps_H/\CQ$. This completes the proof.
\endproof

\begin{example} Consider a homogeneous quantum principal bundle
$P(M,H,\pi)$ with the universal calculus and split by
$i:\ker\eps_H\to \ker\eps_P $.
Let $\Omega^1(H)$ and $\Omega^1(M)$ be determined by $\CQ$ and
$\CN_M$. Then
 there is a
differential calculus $\Omega^1(P)$ with $\Omega^1_{\rm hor}=P(\extd M)P$
and
$\omega(h)=(Si(h)\o)\extd
i(h)\t$
 is a connection on it.
If $\Omega^1(M)$ is left $P$-covariant
then $\Omega^1(P)$ is left-covariant.
The corresponding canonical
map $\ker\eps_H/\CQ\to\ker\eps_P/\CQ_P$ from Proposition~3.5
in this case is $[h]\mapsto \pi_{\CQ_P}\circ i(h)$, where $h\in
\pi^{-1}_{\CQ}([h])\subset\ker\eps_H$. 
\end{example}
\proof We take $\CN_{\rm hor}=\< P\CN_MP,\CN_0\>$  as in Proposition~3.3. Then
$$
\Delta_L(um\tens n v)=u\o (m\tens n)\bo v\o\tens u\t (m\tens n)\bt
v \in P\tens \CN_{\rm hor}
$$ for all $u,v\in P$ and $m\tens n\in
\CN_M$ provided $\Delta_L\CN_M\subset P\tens \CN_M$. Therefore
$\Omega^1_{\rm hor}$ is left-covariant and the assertion follows from
Proposition~3.5. As in Proposition~3.3 the inherited $\Omega^1_P(M)$
is a quotient of $\Omega^1(M)=\Omega^1M/\CN_M$ unless
$\CN_0\cap \Omega^1M\subseteq \CN_M$.\endproof

This provides a natural construction for homogeneous bundles
(where $P$ is a Hopf algebra) to have differential calculi which
are left-covariant.
We conclude with the simplest concrete example of our
construction in Theorem~3.2.

\begin{example} Let $P=H$ regarded as a trivial quantum principal
bundle with $M=\C$ and the universal calculus. The trivialisation is
$\Phi=\id$ and the associated trivial connection is the unique
nonzero $\omega_U$. Hence, for every bicovariant  $\Omega^1(H)$
Theorem~3.2 induces a natural Maurer-Cartan connection
$\omega:\ker\eps/\CQ\to
\Omega^1(H)$.
\end{example}
\proof  Here $\Omega^1M=0$ so $\CN_{\rm hor}=0$ and $\beta=0$ is
the only choice in Proposition~3.3. In fact, there is a unique
connection $\omega_U$ since $\chi=\theta^{-1}$ so that the
condition $\chi\circ\omega_U(h)= 1\tens h$ implies that
$\omega_U(h) =
\theta(1\tens h)$. This is the Maurer-Cartan form with the universal
calculus. We then apply Theorem~3.2. \endproof

\section{Differential Structures on the $q$-Monopole Bundle}

Recall from \cite{BrzMa:gau} that the $q$-monopole (of charge 2) is a
canonical connection in the bundle $SO_q(3)(S_q^2,
\C[Z,Z^{-1}],\pi)$. The quantum group $SO_q(3)$ is a subalgebra of
$SU_q(2)$ spanned by all monomials of even degree. $SU_q(2)$ is
generated by the identity and a matrix
${\bf t} = (t_{ij}) = \pmatrix{\alpha & \beta \cr \gamma &\delta}$, subject
to the homogeneous relations
$$
\alpha\beta = q\beta \alpha , \quad \alpha\gamma = q
\gamma\alpha , \quad \alpha \delta = \delta\alpha +
(q-q^{-1})\beta\gamma ,\quad
\beta\gamma = \gamma\beta, \quad \beta\delta = q\delta \beta
,\quad \gamma \delta = q \delta \gamma ,
$$
and a determinant relation
$\alpha\delta-q\beta\gamma=1$, $q\in \C^*$. We assume that $q$ is not
a root of unity. $SU_ q(2)$ has a matrix quantum
group structure,
$$
\Delta t_{ij} = \sum_{k =1}^ 2 t_{ik}\tens t_{kj}, \quad
\eps(t_{ij})=\delta_{ij},\quad S {\bf t} = \pmatrix{\delta &-q^{-1}
\beta \cr -q\gamma & \alpha}.
$$
The structure quantum group of the $q$-monopole bundle is an algebra
of functions on $U(1)$, i.e. the algebra $\C[Z,Z^{-1}]$ of formal
power series in $Z$ and $Z^{-1}$, where $Z^{-1}$ is an inverse of
$Z$. It has a standard Hopf algebra structure
$$
\Delta Z^{\pm 1} = Z^{\pm 1}\tens Z^{\pm 1} , \quad
\eps(Z^{\pm 1}) = 1, \quad
S Z^{\pm 1} = Z^{\mp 1} .
$$
There is a Hopf algebra projection $\pi: SO_ q(3) \to
k[Z,Z^{-1}]$, built formally from $\pi_{1\over 2} : SU_q(2)\to
\C[Z^{1\over2},Z^{-{1\over 2}}]$,
$$
\pi_{1\over 2} : \pmatrix{\alpha & \beta \cr \gamma &\delta} \mapsto
\pmatrix{Z^{1\over 2} &
0\cr 0 & Z^{-{1\over 2}}},
$$
which defines a right coaction $\Delta_ R : SO_ q(3) \to SO_
q(3) \tens \C[Z,Z^{-1}]$ by $\Delta_ R =(\id\tens\pi)\circ\Delta$.
Finally $S_ q^ 2\subset SO_ q(3)$ is a quantum two-sphere
\cite{Pod:sph}, defined as a fixed point subalgebra, $S_ q^ 2 = SO_
q(3)^{\C[Z,Z^{-1}]}$. $S_ q^ 2$ is generated by $\{1, b_ - =
\alpha\beta , b_ + =\gamma\delta ,b_ 3 = \alpha\delta\}$ and the
algebraic relations in $S_ q^ 2$ may be deduced from those in
$SO_ q(3)$.

The canonical connection in the $q$-monopole bundle $\omega_D$ is
provided  by the
map $i: \C[Z,Z^{-1}] \to SO_q(3)$ given by
$i(Z^{n}) = \alpha^{2n}$,
$i(Z^{-n}) = \delta^{2n}$, $n=0,1,\ldots$ (restricted to
$\ker\eps_{\C[Z,Z^{-1}]}$).
In this section we construct differential structures on the
$q$-monopole bundle using $\omega_D$.

Similarly as in \cite{BrzMa:gau} we choose a differential structure on
$\C[Z,Z^{-1}]$ to be given by
the right ideal $\CQ$ generated by $Z^{-1}+q^4Z-(1+q^4)$. The space
$\ker\eps/\CQ$ is one-dimensional and we denote by $[Z-1]$ its basic
element obtained by projecting $Z-1$ down to $\ker\eps/\CQ$.
\begin{propos}
Let for a quantum principal bundle $SO_q(3)(S_q^2,{
\C}[Z,Z^{-1}],\pi)$, $\CQ$ and $i$ be as above. Then the minimal
horizontal ideal $\CQ_0\in \ker\eps_{SO_q(3)}$ defined in
Proposition~3.5  is generated by the
following elements of $\ker\pi$
$$
\beta\gamma, \qquad q^4\alpha^3\beta +\delta\beta
-(1+q^{4})\alpha\beta , \qquad q^4\alpha^3\gamma +\delta\gamma
-(1+q^{4})\alpha\gamma
$$
The space $\ker\pi/\CQ_0$ and thus the corresponding differential
calculus are infinite-dimensional.

Let $\CQ^{(k,l)}$, $k,l = 1,2,\ldots $ be an infinite family of right
ideals in $\ker\pi$ generated by the generators of $\CQ_0$ and
additionally by 
$\beta^{2k}$, $\gamma^{2l}$. For each pair $(k,l)$,  $\ker\pi/\CQ^{(k,l)}$
is $4(k+l-1)$-dimensional.

Furthermore let $\CQ^{(k,l;r,s)}$,$k,l = 1,2,\ldots $, $r=0,1,\ldots
,k$, $s = 0,1,\ldots ,l$ be an infinite family of right ideals in
$\ker\pi$ generated by the generators of  $\CQ^{(k,l)}$ and also by 
$(\alpha-\delta)\beta^{2r-1}$, $(\alpha-\delta)\gamma^{2s-1}$. Then
$\ker\pi/\CQ^{(k,l;r,s)}$
is a $3k+3l+r+s-4$-dimensional vector space. Notice also that
$\CQ^{(k,l;k,l)} = \CQ^{(k,l)}$.
\end{propos}
\proof The generators of $\CQ_0$ are obtained by a direct computation
of the ideal given in
Proposition~3.5. Explicitly, $\beta\gamma$ is computed by taking
$q=Z^{-1}+q^4Z-(1+q^4)$ and $u=\alpha^2$ and $u=\delta^2$. The
remaining two elements are obtained by taking $q= 1+q^4Z^2-(1+q^4)Z$
and $u=\alpha\beta$ and $u=\alpha\gamma$ correspondingly. It can be then
shown that all the other elements of $\CQ_0$ are generated from the
three listed in the proposition. For example, the choice $q=
Z^{-2}+q^4 -(1+q^4)Z^{-1}$ and $u=\delta\beta$ gives $\delta^3\beta
+q^{4}\alpha\beta
-(1+q^{4})\delta\beta$, but
$$
\delta^3\beta +q^{4}\alpha\beta
-(1+q^{4})\delta\beta = q^{-2}(q^4\alpha^3\beta +\delta\beta
-(1+q^{4})\alpha\beta)\delta^2 - \beta\gamma(q^7\alpha\beta +
q^8\alpha^2\beta\delta -(1+q^4)\beta\delta),
$$
etc. Using this form of the generators of $\CQ_0$ one easily finds
that $\ker\pi/\CQ_0$ is spanned by the projections of the following
elements of $\ker\pi$:
$$
\alpha^k\beta^{2n-k}, \quad \alpha^k\gamma^{2n-k}, \quad
\delta\beta^{2n-1}, \quad \delta\gamma^{2n-1},
\quad n=1,2,\ldots \quad k=0,1,2, \quad k<2n.
$$
Therefore $\ker\pi/\CQ_0$ is an infinite-dimensional vector space.

Notice that for $n=1$ there are 6 independent elements of $\ker\pi/\CQ_0$
coming from monomials in $SO_q(3)$ of degree 1, while for $n>1$ there
are 8 such elements. Using this fact we can compute dimensions of
$\ker\pi/\CQ^{(k,l)}$. Clearly $\dim(\ker\pi/ \CQ^{(1,1)}) =4 = 4(1+1-1)$. Also
$\dim(\ker\pi/\CQ^{(k,l)}) =
\dim(\ker\pi/ \CQ^{(l,k)})$. First take $k=1, l>1$. Then, by counting elements
in $\ker\pi/\CQ_0$ of given degree we find $\dim(\ker\pi/ \CQ^{(1,l)})
= 5 + 4(l-2)+3 =
4l = 4(l+1-1)$. Finally take $k,l>1$. Then $\dim(\ker\pi/\CQ^{(k,l)})=
6+ 4(k-2)
+4(l-2) +6 = 4(k+l-1)$ as stated.

In the case of $\ker\pi/\CQ^{(k,l;r,s)}$ new generators added to
$\CQ^{(k,l)}$ restrict the dimension by $k-r+ l-s$. Therefore
$\dim(\ker\pi/\CQ^{(k,l,r,s)})= \dim(\ker\pi/\CQ^{(k,l)})- (k-r+
l-s) =3k+3l+r+s-4$.
\endproof

\begin{propos}
Let differential structure on $\C[Z,Z^{-1}]$ be given by the ideal
$\CQ$ generated by $Z^{-1}+q^4Z-(1+q^4)$. The largest differential
calculus on $SO_q(3)(S_q^2,{\C}[Z,Z^{-1}],\pi)$ compatible with
q-monopole connection is specified by the ideal $\CQ_P\subset
\ker\eps_{SO_q(3)}$ generated by $\beta\gamma$ and $\delta^2
+q^4\alpha^2 - (1+q^4)$. This calculus is infinite-dimensional.
Let $\CQ_P^{(k,l;r,s)} = \<\CQ^{(k,l;r,s)},i(\CQ)SO_q(3)\>$, be a family of
right ideals in $\ker\eps_{SO_q(3)}$ indexed by $k,l=1,2,\ldots$,
$r=0,1,\ldots , k$, $s=0,1,\ldots ,l$. Each
of $\CQ_P^{(k,l;r,s)}$ induces a $3k+3l+r+s-3$-dimensional, left-covariant
differential calculus on $SO_q(3)$.
\end{propos}
\proof We need to show that $\CQ_P = \<\CQ_0,i(\CQ)SO_q(3)\>$. This
is equivalent to showing that the generators of $\CQ_0$ can be
expressed as linear combinations of elements of $\CQ_P$. Clearly
$\beta\gamma \in \CQ_P$. Furthermore we have
$$
q^4\alpha^3\beta +\delta\beta -(1+q^{4})\alpha\beta = (\delta^2
+q^4\alpha^2 - (1+q^4))\alpha\beta -q^{-3}\beta\gamma\delta\beta \in \CQ_P,
$$
$$
q^4\alpha^3\gamma +\delta\gamma
-(1+q^{4})\alpha\gamma = (\delta^2
+q^4\alpha^2 - (1+q^4))\alpha\gamma -q^{-3}\beta\gamma\delta\gamma
\in \CQ_P.
$$
To prove the remaining part of the proposition it suffices to notice
that $\ker\eps_P/\CQ_P$ is spanned by elements of
$\ker\eps_P/\CQ_0$ listed in Proposition~4.1 and additionally by the
projection of $\alpha^2-1$. Similar calculation as in Proposition~4.1
thus reveals that $\dim(\ker\eps_P/\CQ_P^{(k,l;r,s)}) = 3k+3l+r+s-3$.
\endproof

As a concrete illustration of the above construction we consider differential
calculus induced by $\CQ_P^{(1,1;1,1)}= \CQ_P^{(1,1)}$. Explicitly
$\CQ_P^{(1,1)}$ is
generated by the following four
elements $\delta^2 +q^4\alpha^2 - (1+q^4)$,
$\beta^2,\beta\gamma,\gamma^2$. The space $\ker\eps/\CQ_P^{(1,1)}$ is
five-dimensional, so that $\CQ_P^{(1,1)}$ generates a
five-dimensional left covariant
differential calculus $\Omega^1(SO_q(3))$ on $SO_q(3)$. Since
$SO_q(3)$ is a subalgebra of $SU_q(2)$ the four elements above
generate an ideal in $SU_q(2)$ which also induces a differential calculus
on  $SU_q(2)$.  Choosing the
following basis for the space of left-invariant one forms in
$\Omega^1(SO_q(3))$
\begin{equation}
\omega_0 = {1\over q^4 -1} \pi_{\CN}\circ\theta(1\otimes(q^4\alpha\beta
-\delta\beta)), \qquad \omega_2 = -{q^{-1}\over q^4 -1}
\pi_{\CN}\circ\theta(1\otimes(\delta \gamma
-q^4\alpha\gamma)),
\label{leftforms1}
\end{equation}
\begin{equation}
\omega_3 = {1\over q^2 +1} \pi_{\CN}\circ\theta(1\otimes(\alpha\beta
-\delta\beta)), \qquad \omega_4 = -{1\over q^2 +1}
\pi_{\CN}\circ\theta(1\otimes(\delta\gamma
-\alpha\gamma)),
\label{leftforms2}
\end{equation}
\begin{equation}
\omega_1= {1\over q^{-2} +1}
\pi_{\CN}\circ\theta(1\otimes (\alpha^2-1)),
\label{leftforms3}
\end{equation}
one derives the commutation relations in $\Omega^1(SO_q(3))$ embedded
in $\Omega^1(SU_q(2))$,
\begin{equation}
\omega_{0,2}\alpha = q^{-1}\alpha\omega_{0,2}, \quad \omega_{3,4}\alpha
= q^{-3}\alpha\omega_{3,4}, \quad \omega_1\alpha = q^{-2}\alpha\omega_1
+\beta\omega_4,
\label{rel.omega1}
\end{equation}
\begin{equation}
\omega_{0,2}\beta = q^{1}\beta\omega_{0,2}, \quad \omega_{3,4}\beta
= q^{3}\beta\omega_{3,4}, \quad \omega_1\beta = q^{2}\beta\omega_1
+\alpha\omega_4,
\label{rel.omega2}
\end{equation}
and similarly for $\alpha$ replaced with $\gamma$ and $\beta$ replaced
with $\delta$. The exact one-forms are given in terms of $\omega_i$ as
follows
$$
\extd \alpha = \alpha\omega_1 -q\beta(\omega_2 -{q\over
1-q^2}\omega_4), \qquad \extd\beta =
-q^2\beta\omega_1+\alpha(\omega_0+{q^2\over
1-q^2}\omega_3),
$$
and similarly for $\alpha$ replaced with $\gamma$ and $\beta$ replaced
with $\delta$. It can be easily checked that the forms
$\omega_0,\omega_2,\omega_3,\omega_4$ are horizontal.
Note that this calculus reduces to the 3D calculus of
Woronowicz if one sets $\omega_3 =\omega_4 =0$. This is equivalent to
enlarging $\CQ_P^{(1,1)}$  by $(\delta-\alpha)\beta$,
$(\delta-\alpha)\gamma$ and thus the 3D Woronowicz calculus corresponds to
$\CQ_P^{(1,1;0,0)}$.

The calculus $\CQ_P^{(1,1)}$ appears naturally when one looks at the
monopole bundle from the local point of view. Recall from
\cite{BrzMa:gau} that one of the trivialisations of the $q$-monopole
bundle has the form $P_1(M_1,\C[Z,Z^{-1}],\Phi_1)$,
where
$P_1 = SO_q(3)[(\beta\gamma)^{-1}]$, $M_1 = S_q^2[(b_3-1)^{-1}]$, and
$\Phi_1(Z^n) = (\beta^{-1}\gamma)^n$, $n\in \Z$. This trivialisation
corresponds to the quantum sphere with the north pole removed. It can
be easily shown that $P_1 = M_1\otimes \C[Z,Z^{-1}]$ as an
algebra. The structure of $M_1$ can be most easily described in the
stereographic projection coordinates,
$z=\alpha\gamma^{-1} = qb_-(b_3-1)^{-1}$, $\bar{z} = \delta\beta^{-1}
= b_+(b_3-1)^{-1}$, introduced in \cite{ChuHo:sph}. $M_1$
is then equivalent to the quantum hyperboloid \cite{SchSch:cov} generated
by $z, \bar{z}$, $(1-z\bar{z})^{-1}$ and the relation
$$
\bar{z} z = q^2 z\bar{z}+1 -q^2 .
$$
The natural differential structure $\Omega^1(M_1)$ on $M_1$,
also discussed in \cite{ChuHo:sph}, is given by the relations
$$
z\extd z = q^{-2}\extd z z, \quad z\extd \bar{z} = q^{-2}\extd \bar{z}
z, \quad \bar{z}\extd z = q^2\extd z \bar{z}, \quad \bar{z}\extd
\bar{z} = q^2\extd \bar{z} \bar{z}.
$$
In other words $\Omega^1(M_1) = \Omega^1M_1/ \CN_{M_1}$, where the
subbimodule $\CN_{M_1}\subset \Omega^1M_1$ is generated by
\begin{equation}
q^{-2}\bar{z}\tens z +  z\tens\bar{z} - q^{-2}\bar{z}z\tens 1 -
1\tens z\bar{z}.
\label{zzbar}
\end{equation}
\begin{equation}
(1+ q^2) z\tens z - q^2 z^2\tens 1 - 1\tens z^2, \quad (1+ q^{-2})
\bar{z}\tens \bar{z} - q^{-2}\bar{z}^2\tens 1 - 1\tens \bar{z}^2.
\label{zz}
\end{equation}
The subbimodule $M_1$ and the $q$-monopole connection taken as the
input data in Proposition~\ref{prop.triv} produce the differential
calculus on $P_1$ which coincides with the differential calculus
induced by $\CQ_P^{(1,1)}$ when restricted to $SO_q(3)$. Notice also
that the generator (\ref{zzbar}) appears as a consequence of the
existence of the minimal horizontal subbimodule $\CN_0$. Thus the
differential structures on $P_1$  obtained from data $(\CN_{M_1},\omega_D)$
and $(\tilde{\CN}_{M_1}, \omega_D)$, where $\tilde{\CN}_{M_1}$ is
generated by (\ref{zz}) only, are identical.

In any calculus $\Omega^1(SO_q(3))$ admitting the q-monopole
connection one can define
one-form $\omega_1$ by (\ref{leftforms3}), with $\pi_{\CN}$ a canonical
projection related to the bimodule $\CN$ defining
$\Omega^1(SO_q(3))$. Then the connection
$\omega_D : \ker\eps/\CQ\to\Omega^1(SO_q(3))$ can be computed
explicitly,
$$
\omega_D([Z-1]) = (1+q^{-2})\omega_1.
$$
The canonical map  $i_D: \ker\eps_{\C[Z,Z^{-1}]}/\CQ \to
\ker\eps_{SO_q(3)}/\CQ_P$, with $\CN = \theta(SO_q(3)\otimes\CQ_P)$,
corresponding to $\omega_D$ comes out as
$i_D([Z-1]) = [\alpha^2-1]$ and is clearly $\Ad$-covariant since
$\Ad([Z-1]) = [Z-1]\tens 1$ and $\overline{\Ad}([\alpha^2-1]) =
[\alpha^2-1]\tens 1$.

Similarly, regardless of the differential calculus on
$P_1(M_1,\C[Z,Z^{-1}], \Phi_1)$, the local connection one-form
$\beta:\ker\eps_{\C[Z,Z^{-1}]}\to \Omega^1(M_1)$ can be computed as
follows. It is given by
$$
\beta(h) =
\Phi_1(h\o)Si(h\t)\o\extd(i(h\t)\t\Phi^{-1}_1(h\th)).
$$
To compute it explicitly one can use Lemma~\ref{lemma.cond.beta} to
prove the following
 \begin{lemma}
Let $P(M,\C[Z,Z^{-1}],\Phi)$ be a trivial quantum principal bundle
with a trivialisation $\Phi$ which is an algebra map. Assume that
differential structure $\Omega^1(\C[Z,Z^{-1}])$ is given by the ideal
generated by $Z^{-1}+q^4Z - (1+q^4)$ for $q$ a complex, non-zero
parameter. Then $\omega = \Phi^{-1}*\beta\circ\pi_\eps *\Phi
+\Phi^{-1}*\extd\circ \Phi$ is a connection in
$P(M,\C[Z,Z^{-1}],\Phi)$ if and only if the map $\beta:\ker\eps \to
\Omega^1(M)$ satisfies the following conditions
\begin{eqnarray*}
\beta(Z^{n+1}-1) & = &(1+q^{-4})\Phi(Z)\beta(Z^n-1)\Phi^{-1}(Z) -
q^{-4}\Phi(Z^2)\beta(Z^{n-1}-1)\Phi^{-1}(Z^2) \\
&&+ (1+q^{-4})\Phi(Z)\extd\Phi^{-1}(Z) -
q^{-4}\Phi(Z^2)\extd\Phi^{-1}(Z^2)
\end{eqnarray*}
\begin{eqnarray*}
\beta(Z^{-n}-1) & = &(1+q^{4})\Phi^{-1}(Z)\beta(Z^{-n+1}-1)\Phi(Z) -
q^{4}\Phi^{-1}(Z^2)\beta(Z^{-n+2}-1)\Phi(Z^2)\\
&& +
(1+q^{4})\Phi^{-1}(Z)\extd\Phi^(Z) - q^{4}\Phi^{-1}(Z^2)\extd\Phi(Z^2),
\end{eqnarray*}
for any $n\in \N$.
\end{lemma}

The above lemma implies, in particular, that the map $\beta$
corresponding to the $q$-monopole connection is
fully determined by its action on $Z-1$ say, where it is given by
$$
\beta(Z-1) = (1-z\bar{z})^{-1}(q^2z\extd\bar{z} -
q^{-2}\bar{z}\extd z).
$$
The above formula for $\beta$ is valid in any differential structure
on $M_1$ which admits a $q$-monopole connection, in particular in the natural
one discussed above.
The map $\beta$ is related to $q$-monopole connection $\omega_D$ as
in Proposition~\ref{prop.triv}. The corresponding map $\hat{\Phi}_1$
can be constructed and, applied to the generic element of
$\CN_{\C[Z,Z^{-1}]}$ of
the form $\theta(g\tens h)$, $g\in \C[Z,Z^{-1}]$, $h\in \CQ$, reads
$\Phi_1(g)Si(h)\o\tens i(h)\t$.

\section{Finite gauge theory and Czech cohomology}

In this section we show how quantum differential calculi and gauge
theory can be applied in the simplest setting where $M=\C(\Sigma)$,
$\Sigma$ a finite set, and $H=\C(G)$, $G$ a finite group. We
consider the case of a tensor product bundle
$P=\C(\Sigma)\tens\C(G)$. We show how this formalism provides a
quantum geometrical picture of Czech cohomology when $\Sigma$ is
the indexing set of a good cover of a topological manifold. This
demonstrates a possible new direction to the construction of
manifold invariants: instead of the usual approach in algebraic
topology whereby one looks at the combinatorics of the geometrical
structures on manifolds, we consider instead the (quantum) geometry
of combinatorial structures on the manifold.

Although we are primarily interested in 1-forms (and occasionally
2-forms), it is important to know that they extend to an entire
exterior algebra. Recall that for any unital algebra $M$ there is a
universal extension $\Omega^\cdot M$ of $\Omega^1M$ given in degree
$n$  as the joint kernel in $M^{\tens n+1}$ of all the $n$ maps
given by adjacent product. It can be viewed as
$\Omega^1M\tens_M\Omega^1M\tens_M\cdots\tens_M\Omega^1M$. The collection
$\Omega^\cdot M$ forms a differential graded algebra with
\cmath{(a_0\tens\cdots \tens a_n)\cdot (b_0\tens\cdots \tens b_m)
=(a_0\tens\cdots \tens a_m
b_0\tens\cdots \tens b_m)\\
 \extd_U (a_0\tens\cdots \tens a_n)
=\sum_{j=0}^{n+1}(-1)^{j}(a_0\tens\cdots \tens a_{j-1}\tens 1
\tens a_{j}\tens\cdots \tens a_n)}
with the obvious conventions for $j=0,n+1$ understood. A general
exterior algebra $\Omega^\cdot(M)$ is then obtained by quotienting
it by a differential graded ideal, i.e. an ideal of $\Omega^\cdot
M$ stable under $\extd_U$. Without loss of generality, we always
assume that the degree 0
 component of the differential ideal is trivial. The degree 1 component
 is in particular a sub-bimodule $\CN_M$ of $\Omega^1M$ as in the setting
 above. Conversely,  $\Omega^1(M)$ as defined by a sub-bimodule $\CN_M$ has
 a {\em maximal prolongation} to an exterior algebra $\Omega^\cdot(M)$ by
taking
 differential ideal generated by $\CN_M,\extd_U\CN_M$. In each degree  it can
 be viewed as  a quotient of $\Omega^1(M)\tens_M\Omega^1(M)
 \tens_M\cdots\Omega^1(M)$ by the
additional relations implied by the Leibniz rule applied to the
relations of $\Omega^1(M)$ cf\cite{BrzMa:gau}. For example, $\Omega^2(M)
=\Omega^1(M)\tens_M\Omega^1(M)/ (\pi_M\tens_M\pi_M)(\extd_U\CN_M)$,
where $\pi_M$ is the
canonical projection $\Omega^1M\to\Omega^1(M)$.

Clearly one may take a similar view for $\Omega^2(M)$. The degree 2 part
of a differential ideal of $\Omega^\cdot M$ will, in particular, be a
subbimodule $\CF$ in the range
\[ \overline{\CN_M}\subseteq \CF\subseteq \Omega^2M\]
where $\overline{\CN_M}=(\Omega^1M)\CN_M+\CN_M(\Omega^1M)+\extd_U\CN_M$
is a subbimodule (in view of the Leibniz rule for $\extd_U$), and
$\Omega^2(M)=\Omega^2M/\CF$. Conversely, given $\Omega^1(M)$, any
subbimodule $\CF$ in this range defines an $\Omega^2(M)$ compatible
with $\Omega^1(M)$ in the natural way. Moreover, taking the
differential ideal generated by $\CN_M,\CF,\extd_U\CF$ provides a
prolongation of $\Omega^1(M),\Omega^2(M)$ as specified by
$\CN_M,\CF$. Similarly, one may specify the exterior algebra up to
any finite degree and know that it prolongs to an entire exterior
algebra $\Omega^\cdot(M)$. This is the point of view which we take
throughout the paper.

We begin with a lemma which is well-known (see e.g. \cite[p. 184]{Con:non}),
but which we include because it provides the framework for our
analysis of $\Omega^1$ and $\Omega^2$ in the case of a discrete set.

\begin{lemma} When $\Sigma$ is a finite set of order $|\Sigma|$,
 $\Omega^n\C(\Sigma)$ may be identified with the subset
 $\C^{|\Sigma|}\tens\cdots\tens\C^{|\Sigma|}$ consisting of degree-$(n+1)$
tensors
vanishing on any adjacent diagonal. The exterior derivative
$\Omega^{n-1}\C(\Sigma)\to
\Omega^{n}\C(\Sigma)$ is
\[ (\extd_U
f)_{i_0,\cdots, i_n}=\sum_{j=0}^{n+1}(-1)^{j} f_{i_0,\cdots,
\hat{i_j},\cdots, i_n}\]
where $\hat{\ }$ denotes ommission. The algebra structure of
$\Omega^\cdot\C(\Sigma)$ is $(f\cdot g)_{i_0\cdots
i_{n+m}}=f_{i_0\cdots i_n}g_{i_n\cdots i_{n+m}}$ for $f$ of degree
$n$ and $g$ of degree $m$.
\end{lemma}
\proof  We consider
$\C(\Sigma)$ as a vector space with basis $\Sigma$. An element is
then a vector in $\C^n$ with components $f_i$ for $i\in\Sigma$. The
corresponding function is $f=\sum_i f_i\delta_i$ where $\delta_i$ is
the Kronecker delta-function at $i$. We have $\Omega^n\C(\Sigma)$
as a subspace of $\C(\Sigma)^{\tens n+1}$ in the kernel of adjacent
product maps. These send $\sum f_{i_0\cdots i_n}\delta_{i_0}\tens
\cdots\delta_{i_n}$ to
 $\sum f_{i_0\cdots i_{j-1},i_{j-1},i_{j+1}\cdots i_n}
\delta_{i_0}\tens\cdots\tens\delta_{i_n}$ for all $j=1$ to $j=n$.
So the joint kernel means tensors $f_{i_0,\cdots,i_n}$ vanishing on
the identification of any two adjacent indices. The action of
$\extd_U$ on $\Omega^{n-1}\C(\Sigma)$ is a signed insertion of $1$
in each position of the $n$-fold tensor product, which is the form
stated. The product structure is the pointwise product with the
outer copies of $\C(\Sigma)$, as stated. \endproof

In particular, we identify $\Omega^1\C(\Sigma)$ with $|\Sigma|\times |\Sigma|$
matrices vanishing on the diagonal.

\begin{propos} Let $\Sigma$ be a finite set. Then the possible
$\Omega^1(\C(\Sigma))$
are in 1-1 correspondence with subsets $E\subset
\Sigma\times\Sigma-{\rm Diag}$. The quotient $\Omega^1(\C(\Sigma))$
is obtained by setting to zero the matrix entries $f_{ij}$ for which
$(i,j)\notin E$. In this way we identify $\Omega^1(\C(\Sigma))=\C(E)$.
\end{propos}
\proof
 We consider first the possible sub-bimodules
$\CN_M\subset\Omega^1\C(\Sigma)$. Let $\delta_i$ denote the obvious
(Kronecker delta-function) basis elements of $\C(\Sigma)$. If
$\lambda \delta_i\tens\delta_j+\mu \delta_{i'}\tens\delta_{j'}\in
\CN_M$ for $(i,j)\ne (i',j')$ then multiplying by $\delta_i$ from the
left or by $\delta_j$ from the right implies that $\lambda
\delta_i\tens\delta_j\in
\CN_M$ also, as $\CN_M$ is required to be a sub-bimodule. Hence
$\CN_M={\rm
span}\{\delta_i\tens\delta_j\}$ for $(i,j)$ in some subset of
$\Sigma\times\Sigma-{\rm Diagonal}$. We denote the complement of this subset
in $\Sigma\times\Sigma-{\rm diag}$ by $E$. This gives the general form
of a nonuniversal $\Omega^1(\C(\Sigma))=\Omega^1\C(\Sigma)/\CN_M$.
\endproof

We write $i-j$ whenever $(i,j)\in E$ and we write $i\#j$ whenever
$(i,j)$ is in the complement of $E$ in $\Sigma\times\Sigma-{\rm
diag}$.

\begin{lemma} Let  $\Omega^1(\C(\Sigma))$ be defined as above by $E$.
Then the possible $\Omega^2(\C(\Sigma))$ extending this are in 1-1
correspondence with vector subspaces
$V_{ij}\subset\C(\Sigma-\{i,j\})$ such that
\[ V_{ik}\owns\cases{\sum_{j\ne i,k}\delta_j&if\ $i\#k$\cr
\delta_j&if\ $i\#j, j\ne k$\cr \delta_j& if\ $i\ne j,j\#k$}.\] Then
$\Omega^2(\C(\Sigma))=\oplus_{i,k}\delta_i
\tens{\C(\Sigma-\{i,k\})/V_{ik}}\tens\delta_k$.
We say that $\Omega^2(\C(\Sigma))$ is {\em local} if
all the $V_{ik}$ are spanned by $\delta$-function basis elements.
\end{lemma}
\proof We first compute $\overline{\CN_M}$. Clearly,
 $\CN_M(\Omega^1\C(\Sigma))=
\span\{\delta_i\tens\delta_j\tens\delta_k|\forall i\#j, k\ne j\}$
and  $ (\Omega^1\C(\Sigma))\CN_M=
\span\{\delta_i\tens\delta_j\tens\delta_k|\forall i\ne j, k\# j\}$,
while for $i\#j$,
$\extd_U\delta_i\tens\delta_j=1\tens\delta_i\tens\delta_j-
\delta_i\tens
1\tens\delta_j+\delta_i\tens\delta_j\tens 1$ has most of its terms
contained already in the above. The additional contribution to
$\overline{\CN_M}$ is $\{ \delta_i\tens(\sum_{a\ne
i,j}\delta_a)\tens\delta_j|\ i\#j\}$. These three subspaces span
$\overline{\CN_M}$. Meanwhile, by similar arguments to the proof of
Lemma~5.1, any $\C(\Sigma)$-bimodule $\CF\subset\Omega^2\C(\Sigma)$
has the form
\[ \CF=\span\{\delta_i\tens V_{ik}\tens\delta_k|\ i,k\in\Sigma\},\quad
V_{ik}\subseteq\C(\Sigma-\{i,k\})\] for some vector subspaces as
shown. In order to contain $\overline{\CN_M}$ we see that we require the
subspaces $V_{ik}$ to contain the elements stated. \endproof

The local case is clearly the natural one for `geometry' on the set
$\Sigma$. From Proposition~5.2 we know that $\Omega^1(\C(\Sigma))$
is always local in the same sense. From the above lemma we see that
its maximal prolongation has the $V_{ij}=0$ except in the cases
stated, when it is spanned by the stated vectors; it is therefore
not local and we need to quotient it further.

\begin{theorem} Local $\Omega^2(\C(\Sigma))$ are in correspondence
with subsets
\[ F_0\subseteq\{(i,j)\in\Sigma\times\Sigma|\ i-j,\ j-i\}, \quad F
\subseteq\{(i,j,k)\in
\Sigma\times\Sigma\times\Sigma|\
i-j, j-k, i-k\}.\]
Then $\Omega^2(\C(\Sigma))=\C(F)\oplus\C(F_0)$ can be identified with
3-tensors $f_{ijk}$ vanishing on adjacent diagonals and such that
either $i=k,(i,j)\in F_0$ or $(i,j,k)\in F$.

The 1-cycles in $\Omega^1(\C(\Sigma))$ are $f_{ij}$ such that
\[ f_{ij}=-f_{ji},\quad  f_{ij}-f_{ik}+f_{jk}=0\]
for all $(i,j)\in F_0$ and $(i,j,k)\in F$ respectively. Moreover, the image
of $\C(\Sigma)$ is $(\extd g)_{ij}=g_i-g_j$ for all $i-j$.
\end{theorem}
\proof In the preceding lemma we consider $V_{ik}$ as spanned by
$\delta$-functions on the complement of some subsets
$F_{ik}\subseteq\Sigma-\{i,k\}$, say.  We consider the requirements
of the lemma for the three mutually
exclusive possible cases $i\#k$, $i=k$ and $i-k$. To contain
$\sum_{j\ne i,k}\delta_j$ in the first case, we need
$F_{ik}=\emptyset$.
For the second case, we know that $i\#j$ or $j\#i$ must imply $j$ not in
$F_{ii}$, i.e. $j\in F_{ii}$ should imply $i-j$ and $j-i$ (we consider only
$j\in\Sigma-\{i\}$). This requires $F_{ii}\subset\{j|i-j,j-i\}$.
Similarly for the third possibility. Thus, the conditions on $V_{ik}$ in the
preceding lemma become now
\[ F_{ik}\subseteq \cases{\emptyset&
if\ $i\#k$\cr  \{j\in\Sigma|i-j,j-i\}&if\ $i=k$\cr
\{j\in\Sigma|i-j,j-k\}&if\ $i-k$}.\]
Moreover, in the local case we can identify the quotients as remaining
basis elements, i.e.
$\Omega^2(\C(\Sigma))=\oplus_{i,k}\delta_i\tens\C(F_{ik})\tens\delta_k$.

Next, we can collect together all the $F_{ik}$ where $i-k$. The specification
of
these is equivalent to the specification of $F$ as stated. Likewise,
the specification of all the $F_{ii}$ is equivalent to the specification of
$F_0$ as stated. Then $\Omega^2(\C(\Sigma))=\C(F)\oplus\C(F_0)$ where $\C(F)$
refers to the coefficients of vectors of the form $\delta_i\tens\delta_j
\tens\delta_k$ when $i-j,j-k, i-k$, and $\C(F_0)$ refers to coefficients of
$\delta_i\tens\delta_j\tens\delta_i$.

Finally, we compute the $(\extd
f)_{iji}=f_{ij}+f_{ji}$ and $(\extd f)_{ijk}=f_{ij}-f_{ik}+f_{jk}$
in $\Omega^2(\C(\Sigma))$, where we need only consider $(i,j)\in F_0$ in the
first equation and $(i,j,k)\in F$ in the second. Hence the closed
forms are as stated.
\endproof

It should be clear that a similar situation occurs to all orders. The
maximal prolongation of
local $\Omega^1,\Omega^2$, say, will not be local, requiring further subset
data to obtain local $\Omega^3$, and so on. Note also that such `finite
differential geometry' makes no sense
classically because $1$-forms and functions commute in
$\Omega^1(\C(\Sigma))$ only in the trivial case; one needs the more
general axioms of quantum differential geometry and quantum
exterior algebra.  As an application, we may associate a suitable
nonuniversal quantum differential calculus to any finite cover of a
topological manifold, i.e. we
have the possibility to do `geometry' on the combinatorics of the manifold
rather than combinatorics of the geometry. We
recall that a finite cover $\{U_i\}$ has some nonzero intersections
$\{U_i\cap U_j\}$, some nonzero triple intersections $\{U_i\cap
U_j\cap U_k\}$ etc.

\begin{corol} Let $X$ be a topological manifold with a finite
good open cover $\{U_i\}$ where $i$ run over an indexing set
$\Sigma$. The cover has an associated local quantum
differential calculus $\Omega^1(\C(\Sigma))$, $\Omega^2(\C(\Sigma))$ such
that its quantum cohomology is the Czech cohomology $H^1(X)$.
\end{corol}
\proof   Let $E$ be the distinct pairs for which $U_i\cap
U_j\ne\emptyset$. Here $i-j$ iff $j-i$ so $E$ has a symmetric form. We take
$F_0=E$. We take for $F$ the distinct triples for which
$U_i\cap U_j\cap U_k\ne\emptyset$. We have 1-cochains $\{f_{ij}\}$ defined for
$i-j$ but we do not require $f_{ji}=-f_{ij}$ for the cochain itself, i.e there
are many more 1-cochains than in Czech cohomology. On the other hand, the
closure condition is stronger than in Czech cohomology and antisymmetry appears
`on shell' for any closed cochain. The image of $\extd$ in
$\Omega^1(\C(\Sigma))$
has the usual (antisymmetric) form, so we recover the usual $H^1(X)$ in spite
of
the `quantum' construction.  \endproof

Note that for a smooth compact manifold this recovers the
DeRahm cohomology $H^1(X)$,
i.e. we recover a known geometrical invariant from `geometry'
directly on the cover.
Also, it should be clear that the similar result applies more
generally to any
simplicial complex  (with the one in the corollary being the nerve
of the cover
of a topological manifold.) We let $\Sigma$ be the vertices, $E$ the
edges and $F$ the faces.
The associated
quantum exterior algebra $\Omega^\cdot(\C(\Sigma))$ is such that its
cohomology $H^1$ coincides with the usual simplicial cohomology.
Unlike the usual situation, however, our  `quantum' resolution
of the simplicial cohomology  has the cochains forming a differential graded
algebra and not
only a complex of vector spaces as in the usual situation. This allows us to
proceed in a `geometrical' fashion. Essentially,
the product in Lemma~5.1 is not compatible with antisymmetry of the
cochains and we instead impose the antisymmetry only `on shell'
and not for the cochains themselves.  Although the similarity of $\extd_U$ in
Lemma~5.1
with the Czech coboundary is obvious from the outset, one usually
imposes antisymmetry by hand on the cochains (see for example
\cite{DimMul:non})
and hence loses the exterior algebra structure.

We may now proceed to consider further geometrical structures in this
discrete setting. In
particular, gauge theory or quantum group gauge theory then
provides the natural extension to group or quantum-group valued
Czech cohomology. We note first that if we are interested in only
trivial principal
bundles and gauge theory in terms of the base $M$, we do not need
to fix a differential calculus $\Omega^1(H)$. We need only the coalgebra
structure of
$H$\cite{BrzMa:gau} for a formal gauge theory  with any $\beta:H\to
\Omega^1(M)$ (not
necessarily vanishing on 1) and
any $\gamma:H\to M$ (not necessarily unital). As explained in
\cite{Ma:gosb}
we can use any nonuniversal $\Omega^1(M),\Omega^2(M)$ which are
compatible (as part
of a differential graded algebra), and still have the
fundamental lemma of gauge theory  that
\[ F(\beta)=\extd\beta+\beta * \beta;\quad \beta^\gamma
=\gamma^{-1}*\beta*\gamma+\gamma^{-1}*\extd\gamma\]
obeys
\[ F(\beta^\gamma)=\gamma^{-1}*F(\beta)*\gamma \]
where $*$ denotes the convolution product defined via the coproduct
of $H$. We can still have sections and covariant derivatives as well at
this level\cite{Ma:gosb}. Equally well, we can work with
$\beta\in
\Omega^1(M)\tens A$ and invertible $\gamma\in M\tens A$, where $A$ need only
be a unital algebra. For example, the  zero curvature equation
$\extd \beta+\beta*\beta=0$ makes sense in $\Omega^2(M)\tens A$.

\begin{propos} Let $A$ be a unital algebra and consider gauge fields
$\beta\in \Omega^1(\C(\Sigma))\tens A$  such that $F(\beta)=0$ in
$\Omega^2(\C(\Sigma))\tens A$.
There is an action of the group of invertible elements $\gamma\in
\C(\Sigma)\tens A$ on this space and the moduli space of zero curvature
gauge fields modulo such transformations coincides with
the multiplicative Czech cohomology $H^1(X,A)$ in the setting of the preceding
corollary.
\end{propos}
\proof In the setting of Proposition~5.4 we have
\[ F(\beta)_{iji}=\beta_{ij}+\beta_{ji}+\beta_{ij}\beta_{ji},\quad
F(\beta)_{ijk}
=\beta_{ij}+\beta_{jk}-\beta_{ik}
+\beta_{ij}\beta_{jk}\]
for all $(i,j)\in F_0$ and $(i,j,k)\in F$ respectively. Hence the
zero-curvature equation is
\[ (1+\beta_{ij})(1+\beta_{ji})=1, \quad (1+\beta_{ij})(1+\beta_{jk})
=1+\beta_{ik}\]
as a multiplicative version of Proposition~5.4 and with values in
$A$. Although $\beta_{ij}$ are not imposed to be such that
$g_{ij}=1+\beta_{ij}$
is invertible, we see that this appears `on shell' for zero curvature gauge
fields,
along with $g_{ij}^{-1}=g_{ji}$. Finally, a gauge transformation means
$\gamma\in \C(\Sigma)\tens A$ with components $\{\gamma_i\}$ invertible, and
the action on
connections is
\[ \beta^\gamma_{ij}=\gamma_i^{-1}\beta_{ij}\gamma_j
+\gamma_i^{-1}\gamma_j-1\]
for all $i-j$. Hence, in the particular setting of Corollary~5.5 (or more
generally for a simplicial complex) we obtain for the moduli space of zero
curvature gauge fields the multiplicative Czech cohomology.   \endproof

Note that if $A$ supports logarithms then $1+\beta_{ij}=\exp{f_{ij}}$ and the
multiplicative theory becomes equivalent to the additive theory as in
Corollary~5.5, i.e. we have a second interpretation with $f$ as $A$-valued
quantum differential forms in this case.

We proceed now to quantum group gauge theory with a full quantum
geometric structure where $P=\C(\Sigma)\tens
\C(G)=\C(\Sigma\times G)$, $G$ a finite group (say) and both
$\C(G)$, $\C(\Sigma)$ are equipped with quantum differential
calculi. Bicovariant (coirreducible) calculi on $\C(G)$ are known to
correspond to nontrivial conjugacy classes on $G$. When $G=\Z_2$
there is only one non-zero calculus, which is also the universal
one. Here $\ker\eps$ is 1-dimensional so $\beta,\gamma$ are fully
specified as $\beta\in\Omega^1(\C(\Sigma))$ and
$\gamma\in\C(\Sigma)$ with invertible components. In this case we
recover the setting of Proposition~5.6 with $A=\C$. However, for
other groups (or if we use the zero calculus on $\C(\Z_2)$) we need
the theory of quantum principal bundles with nonuniversal calculi
developed in Section~3. We demonstrate some of this theory now,
namely Proposition~3.3 which provides the construction of the
differential calculus on a trivial bundle $P$ by `gluing' the
chosen calculi on the base and on the fibre via a universal
connection.

We consider $G=\Z_3=\{e,g,g^2\}$, which has two non-zero
bicovariant calculi, associated to $g$ or $g^{-1}$. Without loss of
generality we consider the one associated to $g$. Then
$\Omega^1(\C(\Z_3))$ is 1-dimensional over $H=\C(\Z_3)$. The unique
normalised left-invariant 1-form is $\omega_1$ say and
\[\extd\delta_e=(\delta_{g^2}-\delta_e)\omega_1,\quad \extd\delta_g
=(\delta_e-\delta_g)\omega_1,
\quad \omega_1 \delta_{g^i}=\delta_{g^{i-1}}\omega_1\]
gives its structure on a $\delta$-function basis of $\C(\Z_3)$. The
ideal $\CQ$ for this bicovariant calculus is $\CQ=\C\delta_{g^2}$.
{}From the point of view of Proposition~5.2, the calculus corresponds
to edges specified by $a-b$ iff $a=b-1$, where $a,b\in\{0,1,2\}$
mod 3. The corresponding subbimodule of $\Omega^1\C(\Z_3)$ is
$\span\{\delta_e\tens\delta_{g^2},\delta_g\tens\delta_e,
\delta_{g^2}\tens\delta_g\}$.

\begin{example} Let $\C(\Sigma)$ have differential calculus described by
a collection of edges $\{i-j\}$ via Proposition~5.2. Let $\C(\Z_3)$
have the standard 1-dimensional calculus as above. For any
$\beta_U:\ker\eps\to\Omega^1\C(\Sigma)$, i.e. a pair
$\beta^{(1)}=\beta_U(\delta_{g})$,
$\beta^{(2)}=\beta_U(\delta_{g^2})$ of $|\Sigma|\times|\Sigma|$ of
matrices with zero diagonal, the induced $\Omega^1(\C(\Sigma\times
\Z_3))$ via Proposition~3.3 has the allowed
edges
\align{(i,a)-(j,a)&{\rm if}& i-j,\ \beta^{(2)}_{ij}=0\\
             (i,a)-(i,b)&{\rm if}&   a=b-1\\
             (i,a-1)-(j,a)&{\rm if}& i-j,\ \beta^{(1)}_{ij}=0}
Moreover, $\omega:\ker\eps/\CQ\to\Omega^1(\C(\Sigma\times \Z_3))$
defined by
\align{\omega(\delta_{g})\equad&&=\sum_{i-j,\beta^{(2)}_{ij}=0}
\sum_a \beta^{(1)}_{ij}\delta_i\tens\delta_{g^a}
\tens\delta_j\tens\delta_{g^a}
-\sum_{i-j,\beta^{(1)}_{ij}=0}\sum_a\delta_i\tens
\delta_{g^{a-1}}\tens\delta_{j}
\tens\delta_{g^a}\\&&-\sum_{i,a}\delta_i\tens\delta_{g^{a-1}}
\tens\delta_i\tens\delta_{g^a} }
is a connection on $\C(\Sigma\times \Z_3)$ as a quantum principal
bundle with this quantum differential calculus.
\end{example}
\proof Since $P=\C(\Sigma)\tens\C(\Z_3)$ is a tensor product
bundle $P=M\tens H$, the
trivialisation in Proposition~3.3 is
$\Phi(h)=1\tens h$ and so
\[ \omega_U(h)=(1\tens Sh\o)\beta_U(\pi_\eps(h\t))\tens h\th+1\tens
Sh\o\tens 1\tens h\t 
-1\tens1\tens1\tens 1\eps(h).\]
To compute the minimal horizontal subbimodule
\[\CN_0=P\span\{(m\tens h\o)\omega_U(qh\t)-\omega_U(q)(m\tens h)
|\ q\in\CQ, \ m\in M,\ h\in H\}\]
and $\CN=\<P\CN_MP,P\omega_U(\CQ)P\>$ defining $\Omega^1(P)$, we
compute first
\align{\omega_U(\delta_{g^2})\equad&&=
\sum_{a+b+c=2}\delta_{g^{-a}}\beta_U(\delta_{g^b})\delta_{g^c}
+\sum_{a+b=2}1\tens\delta_{g^{-a}}\tens 1\tens\delta_{g^b}\\
&&=\sum_{i,j,a}\beta^{(1)}_{ij}\delta_i\tens\delta_{g^{a-1}}
\tens\delta_j\tens\delta_{g^{a}}
+\sum_{i,j,a}\beta^{(2)}_{ij}\delta_i\tens\delta_{g^a}\tens
\delta_j\tens\delta_{g^a}
-\sum_{a}1\tens\delta_{g^{a+1}}\tens1\tens\delta_{g^{a}}}
where indices $a,b,c$ are taken in $\{0,1,2\}$ mod 3 and
$i,j\in\Sigma$. Then
\align{\equad&& (\delta_l\tens\delta_{g^b})\omega_U(\delta_{g^2})
(\delta_k\tens\delta_{g^{a}})\\
&&=(\delta_l\tens\delta_{g^b}\tens1\tens 1)\times\\&&\times(\sum_{j}
\beta^{(1)}_{jk}\tens\delta_j\tens\delta_{g^{a-1}}
\tens\delta_k\tens\delta_{g^a}
+\sum_j\beta^{(2)}_{jk}\delta_j\tens\delta_{g^a}\tens
\delta_k\tens\delta_{g^a}
-1\tens\delta_{g^{a+1}}\tens\delta_k\tens\delta_{g^a})\\
&&=\delta_{b,a-1}\beta^{(1)}_{lk}\delta_l\tens
\delta_{g^{a-1}}\tens\delta_k\tens\delta_{g^a}
+\delta_{b,a}\beta^{(2)}_{lk}\delta_l\tens\delta_{g^a}
\tens\delta_k\tens\delta_{g^a}
+\delta_{b,a+1}\delta_l\tens\delta_{g^{a+1}}\tens\delta_k\tens\delta_{g^a}}
Choosing $b=a-1,a,a+1$ we see that
\align{P\omega_U(\CQ)P&&\equad
=\span\{\delta_{i}\tens\delta_{g^{a-1}}
\tens\delta_j\tens\delta_{g^a}|\beta^{(1)}_{ij}\ne0\}
+\span\{\delta_{i}\tens\delta_{g^{a}}\tens\delta_j\tens\delta_{g^a}|
\beta^{(2)}_{ij}\ne0\}\\
&&\quad+\span\{\delta_{i}\tens\delta_{g^{a+1}}
\tens\delta_j\tens\delta_{g^a}\}.}
This and
\[ P\CN_MP=\span\{\delta_{i}\tens\delta_{g^{a}}
\tens\delta_j\tens\delta_{g^b}|\ i\#j\}\]
gives $\CN$. One may compute $\CN_0$ similarly, noting that since
$\CQ=\C\delta_{g^2}$,
\[\CN_0=\span\{(\delta_l\tens\delta_{g^{b}}\tens1\tens 1)\left((\delta_k\tens
\delta_{g^{a+1}})\omega_U(\delta_{g^2})
-\omega_U(\delta_{g^2})(\delta_k\tens\delta_{g^a})\right)\}.\]
This turns out to be the $P\omega_U(\CQ)P$ in which its third part
is restricted to
$\span\{\delta_i\tens\delta_{g^{a+1}}\tens\delta_j\tens\delta_{g^a}|i\ne
j\}$.

Next, we compute the edges corresponding to $\CN$ as in the setting
of Proposition~5.2. We  consider only $(i,a)\ne (j,b)$. Then
$(i,a)\#(j,b)$ whenever $a=b+1$ or $(a=b-1,\ \beta^{(1)}_{ij}\ne
0)$ or $(a=b,\ \beta^{(2)}_{ij}\ne 0)$. So the complementary set is
$(i,a)-(j,b)$ whenever $(a=b\ {\rm or}\ a=b-1)$ and $(i=j\ {\rm
or}\ i-j)$ and $(a=b\ {\rm or}\
\beta^{(1)}_{ij}=0)$ and $(a=b-1\ {\rm or}\ \beta^{(2)}_{ij}=0)$,
which simplifies as stated.

Finally, Proposition~3.3 also provides for a connection
$\omega:\ker\eps/\CQ\to\Omega^1(P)$. In our case we identify
$\ker\eps/\CQ=\C\delta_{g}$. Then
\align{ \omega_U(\delta_g)
\equad &&=\sum_{i,j,a}\beta^{(1)}_{ij}\delta_i\tens
\delta_{g^a}\tens\delta_j\tens\delta_{g^a}
+\sum_{i,j,a}\beta^{(2)}_{ij}\delta_i\tens\delta_{g^{a+1}}
\tens\delta_j\tens\delta_{g^a}-\sum_a
1\tens \delta_{g^{a-1}}\tens 1\tens \delta_{g^a}.} We then project
this down by setting to zero elements in $\CN$, which gives the
result as shown. In specific examples one may also compute
$\Omega^1_P(M)$ obtained by restricting $\Omega^1(P)$ to $M$ (in
general it will not be our original $\Omega^1(M)$, having instead
the new subbimodule $\Omega^1M\cap\CN$).
\endproof

We see that a connection $\beta_U$ `glues' the differential
calculus in $\C(G)$ to that on $\C(\Sigma)$ to obtain a
differential calculus on the total space. We can of course take
quantum groups other than $\C(G)$. For example, we may take $H=\C
G$, $G$ a finite group. When $G$ is non Abelian, $H$ is not the
function algebra on any space, so this is a genuine application of
`noncommutative geometry'. In this case we know from \cite{Ma:cla}
that (coirreducible) bicovariant calculi $\Omega^1(\C G)$ may be
identified with pairs $(V,\lambda)$ where $V$ is an (irreducible)
representation and $\lambda\in P(V^*)$.  We will construct
nonuniversal calculi and connections on bicrossproduct bundles of
this type (i.e. with fibre $\C G$) in the next section.

One can (in principle) consider other connections on this bundle,
the zero curvature condition etc., and obtain in this way (in view
of Proposition~5.6) a slew of refinements of Czech cohomology with
values in quantum groups equipped with quantum differential
structures. Recall that at the level of naive gauge theory as in
Proposition~5.6 only the coalgebra of $H$ enters. Thus $H=\C G$
just yields $|G|-1$ copies of the 1-dimensional gauge theory. By
contrast, the theory with nonuniversal calculi on the fibre and
bundle carries much more information, including the group structure
and (in the case of $\C G$) the choice of $(V,\lambda)$. One also
has extensions of the geometric theory of quantum principal bundles
where the fibre is a braided group or only a
coalgebra\cite{BrzMa:coa}\cite{Ma:diag}. In a dual form it means
gauge fields with values in algebras (not necessarily Hopf
algebras) equipped with differential calculi.

Finally, the extension of these ideas to values in a sheaf is also
important. Valuation of the usual Czech $H^1$ in the structure
sheaf provides of course a classification of line bundles over $X$,
etc. By taking more exotic Hopf algebras and differential calculi
in a sheaf setting we may obtain more interesting invariants and
`quantum geometrical' methods to compute them. A further long-range
suggestion provided  by the above result is that the role of an
`open cover' can be naturally encoded as a discrete algebra (here
$\C(\Sigma)$) and the choice of nonuniversal differential calculus
on it. One may be able to turn this around and take a discrete
algebra $M$ and choice of $\Omega^\cdot(M)$ on it as the starting
point for the definition of a quantum manifold `with cover $M$'.
One should then define a `sheaf over $M,\Omega^\cdot(M)$', etc.
These are directions to be explored elsewhere.

\section{Bundles and Connections on Cross Product Hopf Algebras}

As noted already in Section~2, a general trivial quantum principal
bundle has the form of a cocycle cross product. Here we will
consider in detail some special cases of such cross products where
the total space $P$ is itself a Hopf algebra. This covers many of
the Hopf algebras in the literature, providing for them natural
calculi and connections. This is a further concrete application of
quantum group gauge theory and provides a uniform approach to the
different kinds of cross product.

In fact, there are mainly two different general constructions for
Hopf algebras where the algebra part is a cross product. The first,
the bicrossproduct construction\cite{Ma:phy} associates quantum
groups to group factorisations. The other is a bosonisation
construction\cite{Ma:bos} which provides the Borel and maximal
parabolic parts of the quantum groups $U_q(\cg)$, as well as a way
of thinking about the quantum double\cite{Ma:skl}\cite{Ma:mec} and
Poincar\'e quantum groups\cite{Ma:poi}. Slightly more general is a
biproduct construction\cite{Rad:str}\cite{Ma:skl}, with the
starting point being a braided group.

Note that if a homogeneous bundle as in Example~2.3 is split by a
coalgebra map $i:H\to P$ then (a) the bundle is trivial by $\Phi=i$
and (b) the $\Ad$-invariance condition in Example~2.3 holds and the
canonical connection $Si(h)\o\extd i(h)\t$ coincides with the
trivial $\beta=0$ connection in Example~3.3. The bosonisations are
of this type (in fact, $i$ a Hopf algebra map), while
bicrossproducts are not in general of this type, although the
bundle is still trivial.

\subsection{Bicrossproducts}

We recall \cite{Ma:book} that a general extension of Hopf algebras
has the form of a bicrossproduct
\[ M\to M\bicross H\to H\]
possibly with cocycles.

We consider  the   cocycle-free case. In this case $H$ acts on $M$
and $M$ coacts on $H$ and the Hopf algebra structure is the
associated cross product and cross coproduct (or `bicrossproduct')
from \cite{Ma:phy}. In this case it is immediate to see from the
explicit formulae that $\pi: M\bicross H\to H$, $\pi(m\tens
h)=\eps(m)h$ is a homogeneous quantum principal bundle over $M$.
Moreover, the map $\Phi:H\to M\bicross H$, $\Phi(h)=1\tens h$ is an
algebra map. It is easy to see that
$\Delta_R\circ\Phi=(\Phi\tens\id)\circ\Delta$ and
$\Phi^{-1}=\Phi\circ S$, so that $M\bicross H$ as a bundle is
trivial. From Proposition~3.3 we already know that natural calculi
$\Omega^1(P)$  are provided by the choice of connection defined by
$\beta_U:\ker\eps\to\Omega^1M$. We provide now a construction for
suitable $\beta_U$ such that the resulting $\Omega^1(P)$ is
left-invariant.

\begin{propos} Strong, left-invariant
connections in $M\bicross H$ as
a trivial quantum principal bundle are in 1-1 correspondence with
linear left-invariant maps $\beta_U:\ker\eps\to\Omega^1M$ such that
\[ \beta_U(\pi_\eps(h\o))h\t\bt\tens h\t\bo - h\o\bt
\beta_U(\pi_\eps(h\t))\tens h\o\bo =  \extd_Uh\bt\tens h\bo\]
Moreover,
such $\beta_U$ are in 1-1 correspondence with linear maps
$\gamma:H\to M$ obeying $\gamma(1)=1$ and $\eps_M\circ\gamma =
\eps_H$, and such that
\[ \gamma(h\o)h\t\bt\tens h\t\bo=\gamma(h\t)\tens h\o, \quad\forall h\in H.\]
The correspondence is via
\[ \beta_U(h)=(S\gamma(h)\o)\extd_U \gamma(h)\t.\]
The corresponding $\omega_U$ is a canonical connection for a splitting map
$i(h)=\gamma(h\o)\tens h\t$.
\end{propos}
\proof First recall some basic facts about  bicrossproducts that are
needed for the proof. The definition of a coproduct in $M\bicross H$,
$\Delta(m\tens h) = m\o\tens h\o\bo\tens m\t h\o\bt \tens h\t$ implies
that $\Delta\Phi(h) = \Phi(h\o\bo)\tens h\o\bt\Phi(h\t)$. Let $\alpha
: H\to M\tens H$ denote a right coaction of $M$ on $H$ used for the
definition of $M\bicross H$, i.e. $\alpha(h) = h\bo\tens h\bt$. Then
the property $\alpha(gh) = g\o\bo h\bo\tens g\o\bt(g\t\la h\bt)$
implies
\begin{eqnarray*}
1_H\tens \Phi(Sh) &=& \alpha(Sh\t h\th)\Phi(Sh\o) = Sh\th\bo h\fo\bo
\tens Sh\th\bt(Sh\t\la h\fo\bt)\Phi(Sh\o) \\
& = & Sh\fo\bo h\fiv\bo\tens Sh\fo\bt\Phi(Sh\th)h\fiv\bt\Phi(S^2h\t
Sh\o) \\
& = & Sh\t\bo h\th \bo \tens Sh\t\bt\Phi(Sh\o)h\th\bt,
\end{eqnarray*}
where we used the fact that $h\la m = \Phi(h\o)m\Phi(Sh\t)$,
$\forall m\in M, h\in H$. Therefore
\begin{equation}
1_H\tens \Phi(Sh)
=  Sh\t\bo h\th \bo \tens Sh\t\bt\Phi(Sh\o)h\th\bt .
\label{bicross.2}
\end{equation}
Similarly, for any $h\in H$
\begin{equation}
\eps(h)1_H\tens 1_M = h\o\bo Sh\fo\bo\tens
h\o\bt\Phi(h\t)Sh\fo\bt\Phi(Sh\th).
\label{bicross.1}
\end{equation}
Now we can start proving the proposition. First
assume that $\omega_U$ is a strong, left-invariant connection. Recall
from \cite{Haj:str} that the connection $\Pi :\Omega^1P\to\Omega^1P$
in $P(M,H)$
is said to be strong if $(\id-\Pi)(\extd_U P) \subset (\Omega^1 M)P$. In
the case of
a trivial bundle $P(M,H,\Phi)$ this is equivalent to the existence of
a map $\beta_U :\ker\eps_H \to \Omega^1M$, given by  $\beta_U(h) =
\Phi(h\o)\omega_U(\pi_\eps(h\t))\Phi^{-1}(h\th) +
\Phi(h\o)\extd_U\Phi^{-1}(h\t)$. In our case $\Phi$ is an algebra map,
therefore $\Phi^{-1} = \Phi\circ S$. Since $\omega_U$ is assumed to be
left-invariant we find, for any $h\in\ker\eps_H$,
\begin{eqnarray*}
\Delta_L\!\!\!\!\!\!\!\!&
&\!\!\!\!\!\!\!\!(\beta_U(\pi_\eps(h\o))h\t\bt)\tens h\t\bo\\
&  = &
\Phi(h\o)\o\Phi(Sh\th)\o h\fo\bt\o\tens
\Phi(h\o)\t\omega_U(\pi_\eps(h\t)) \Phi(Sh\th)\t h\fo\bt\t\tens h\fo\bo\\
&& + \Phi(h\o)\o\Phi(Sh\t)\o h\th\bt\o\tens
\Phi(h\o)\t(\extd_U\Phi(Sh\t)\t)h\th\bt\t\tens h\th\bo \\
& = & \Phi(h\o\bo Sh\fiv\bo)h\six\bt\o\tens
h\o\bt\Phi(h\t)\omega_U(\pi_\eps(h\th))
Sh\fiv\bt\Phi(Sh\fo)h\six\bt\t\tens h\six\bo \\
&& +\Phi(h\o\bo Sh\fo\bo)h\fiv\bt\o\tens
h\o\bt\Phi(h\t)\extd_U(Sh\fo\bt\Phi(Sh\th))h\fiv\bt\t\tens h\fiv\bo.
\end{eqnarray*}
On the other hand since $\beta_U(h)\in \Omega^1 M$,
$\Delta_L(\beta_U(h))\in M\tens \Omega^1M$, i.e. $\Delta_L$ is the
coaction of $M$ on $\Omega^1M$. Therefore the outcome of
the above calculation must be in $M\tens \Omega^1M\tens H$. Applying
$1_H\eps_M\tens \id_{\Omega^1M}\tens \id_H$ and noting that
$(\eps_M\tens \id)(1\tens h)(m\tens 1)
= \eps(m)h$, for any $h\in H$ and $m\in M$ we find
\begin{eqnarray*}
1_H\tens  \beta_U(\pi_\eps(h\o))h\t\bt\!\!\!\!\!\!\!\!&
&\!\!\!\!\!\!\!\!\tens h\t\bo\\
 & = & h\o\bo
Sh\fiv\bo\tens
h\o\bt\Phi(h\t)\omega_U (\pi_\eps(h\th))Sh\fiv\bt\Phi(Sh\fo)h\six\bt\tens
h\six\bo \\
&& +h\o\bo Sh\fo\bo\tens
h\o\bt\Phi(h\t)\extd_U(Sh\fo\bt\Phi(Sh\th)h\fiv\bt)\tens h\fiv\bo \\
&& -h\o\bo Sh\fo\bo\tens
h\o\bt\Phi(h\t)Sh\fo\bt\Phi(Sh\th)\extd_Uh\fiv\bt\tens h\fiv\bo .
\end{eqnarray*}
This implies
\begin{eqnarray*}
\beta_U\!\!\!\!\!\!\!\!&
&\!\!\!\!\!\!\!\!(\pi_\eps(h\o))h\t\bt\tens h\t\bo\\
 & = & h\o\bt\Phi(h\t)\omega_U
(\pi_\eps(h\th))Sh\fiv\bt\Phi(Sh\fo)h\six\bt\tens
h\o\bo Sh\fiv\bo h\six\bo \\
&& + h\o\bt\Phi(h\t)\extd_U(Sh\fo\bt\Phi(Sh\th)h\fiv\bt)\tens h\o\bo
Sh\fo\bo h\fiv\bo \\
&& -h\o\bt\Phi(h\t)Sh\fo\bt\Phi(Sh\th)\extd_Uh\fiv\bt\tens h\o\bo
Sh\fo\bo h\fiv\bo  \\
& = & h\o\bt\Phi(h\t)\omega_U(\pi_\eps(h\th))\Phi(Sh\fo)\tens h\o\bo +
h\o\bt\Phi(h\t)\extd_U\Phi(Sh\th)\tens h\o\bo \\
&& +h\o\bt\eps(h\t)\tens 1 \tens h\o\bo -
h\o\bt\Phi(h\t)Sh\fo\bt\Phi(Sh\th)\tens h\fiv\bt\tens h\o\bo
Sh\fo\bo h\fiv\bo \\
& = & h\o\bt\beta_U(\pi_\eps(h\t))\tens h\o\bo + h\bt\tens 1\tens h\bo
-\eps(h\o)\tens h\t\bt\tens h\t\bo \\
& = & h\o\bt
\beta_U(\pi_\eps(h\t))\tens h\o\bo - \extd_Uh\bt\tens h\bo,
\end{eqnarray*}
where we used property (\ref{bicross.2}) and definition of the
universal differential to derive the second equality and
(\ref{bicross.1})  to derive the third one.

Furthermore, we find
\begin{eqnarray*}
\Delta_L(\beta_U(h)) &= & \Phi(h\o)\o\Phi(Sh\th)\o\tens
\Phi(h\o)\t\omega_U(\pi_\eps(h\t)) \Phi(Sh\th)\t\\
&& + \Phi(h\o)\o\Phi(Sh\t)\o\tens
\Phi(h\o)\t\extd_U\Phi(Sh\t)\t \\
& = & \Phi(h\o\bo Sh\fiv\bo)\tens
h\o\bt\Phi(h\t)\omega_U(\pi_\eps(h\th)) Sh\fiv\bt\Phi(Sh\fo) \\
&& +\Phi(h\o\bo Sh\fo\bo)\tens
h\o\bt\Phi(h\t)\extd_U(Sh\fo\bt\Phi(Sh\th))
\end{eqnarray*}
Using the fact that
$\Delta_L(\beta_U(h))\in M\tens \Omega^1M$ and that $M$ is invariant
under $\Delta_R$ we can apply $\Delta_R$ to first factor in
$\Delta_L(\beta_U(h))$ then $\Phi^{-1}$ to second factor in the
resulting tensor product and multiply first two factors to obtain back
$\Delta_L(\beta_U(h))$. Applying the same procedure to the right hand
side of the above equality, using the fact that $\Phi$ is an
intertwiner for the right coaction of $H$ on $M\bicross H$ as well as
the properties of a counit in $M\bicross H$ we thus find
\begin{eqnarray*}
\Delta_L(\beta_U(h)) &= & \eps(h\o\bo Sh\fiv\bo)1\tens
h\o\bt\Phi(h\t)\omega_U(\pi_\eps(h\th)) Sh\fiv\bt\Phi(Sh\fo) \\
&& +\eps(h\o\bo Sh\fo\bo)1\tens
h\o\bt\Phi(h\t)\extd_U(Sh\fo\bt\Phi(Sh\th))\\
& = & 1\tens (\Phi(h\o)\omega_U(\pi_\eps(h\t))\Phi(Sh\th) +
\Phi(h\o)\extd_U\Phi(Sh\t)) = 1\tens\beta_U(h).
\end{eqnarray*}
Therefore $\beta_U$ is left-invariant as stated.

Conversely, let $\beta_U:\ker\eps_H\to \Omega^1 M$ be a left-invariant
linear map satisfying the condition in the proposition. Define
$\omega_U :\ker\eps_H \to \Omega^1P$ by 
$\omega_U(h) = \Phi(Sh\o)\beta_U(\pi_\eps(h\t))\Phi(h\th) +
\Phi(Sh\o)\extd_U\Phi(h\t)$. The map $\omega_U$ is a strong connection
1-form. We need to verify whether it is left-invariant. For any
$h\in\ker\eps_H$ 
we use the left-invariance of $\beta_U$ and compute
\begin{eqnarray*}
\Delta_L\omega(h) &=& \Phi(Sh\t\bo h\fo\bo)\tens Sh\t\bt
\Phi(Sh\o)\beta_U(\pi_\eps(h\th))h\fo\bt\Phi(h\fiv)\\
&& + \Phi(Sh\t\bo h\th\bo)\tens Sh\t\bt
\Phi(Sh\o)\extd_U(h\th\bt\Phi(h\fo))\\
&=& \Phi(Sh\t\bo h\th\bo)\tens Sh\t\bt
\Phi(Sh\o)h\th\bt\beta_U(\pi_\eps(h\fo))\Phi(h\fiv)\\
&& - \Phi(Sh\t\bo h\th\bo)\tens Sh\t\bt
\Phi(Sh\o)(\extd_U h\th\bt)\Phi(h\fo)\\
&& + \Phi(Sh\t\bo h\th\bo)\tens Sh\t\bt
\Phi(Sh\o)(\extd_Uh\th\bt)\Phi(h\fo)\\
&& + \Phi(Sh\t\bo h\th\bo)\tens Sh\t\bt
\Phi(Sh\o)h\th\bt\extd_U\Phi(h\fo)\\
&=& 1\tens (\Phi(Sh\o)\beta_U(\pi_\eps(h\t))\Phi(h\th)+
\Phi(Sh\o)\extd_U\Phi(Sh\t)),
\end{eqnarray*}
where the assumption about $\beta_U$ and the Leibniz rule were used in
 the derivation of the 
 second equality and the property (\ref{bicross.2}) in derivation of
 the last one. Therefore $\omega_U$ is a left-invariant connection as
 required.

Since $\beta_U(h)$ is a left-invariant form on $M$ for any
$h\in\ker\eps_H$ then the similar argument as in the proof of
Proposition~3.4 yields that $\beta_U(h) = S\gamma(h)\o\extd_U
\gamma(h)\t$ with $\gamma = (\eps_M\tens \id)\circ\beta_U$, a map
$\ker\eps_H\to\ker\eps_M$,  which is
extended uniquely to $H$ by setting $\gamma(1) = 1$. In other words
$\gamma(h) = (\eps_M\tens \id)\circ\beta_U(\pi_\eps(h)) +\eps(h)1_M$,
for any $h\in H$. Notice that $\eps_M(\gamma(h)) = \eps_H(h)$.  Assuming
that $\beta_U$ satisfies the condition specified in the proposition
and applying $\eps_M\tens \id$ one finds
$$
(\gamma(\pi_\eps(h\o))+\eps(h\o))h\t\bt\tens
h\t\bo=(\gamma(\pi_\eps(h\t))+\eps(h\t))\tens h\o,
$$
i.e.
\[ \gamma(h\o)h\t\bt\tens h\t\bo=\gamma(h\t)\tens h\o, \quad\forall
h\in H,\]
as required. Now take any map $\gamma :H\to M$, $\gamma(1)=1$,
$\eps_M\circ\gamma = \eps_H$, and  such
that the above condition is satisfied. Applying
$(S\tens\id)\circ\Delta$ to the first factor in this equality and
using definition of the universal differential one finds
\begin{eqnarray*}
&&Sh\t\bt\o S\gamma(h\o)\o(\extd_U\gamma(h\o)\t)h\t\bt\t\tens h\t\bo -
 S\gamma(h\t)\o\extd_U\gamma(h\t)\t\tens h\o\\
&&  =
 Sh\bt\o\extd_Uh\bt\t\tens h\bo,
\end{eqnarray*}
or, by using the form of $\beta_U$, i.e. $\beta_U(h)
= S\gamma(h)\o\extd_U\gamma(h)\t$
$$
Sh\t\bt\o \beta_U(\pi_\eps(h\o))h\t\bt\t\tens h\t\bo -
 \beta(\pi_\eps(h\t))\tens h\o =
 Sh\bt\o\extd_Uh\bt\t\tens h\bo.
$$
By applying the coaction $\alpha$ to the second factor in the above
equality, interchanging
third factor with the second and the first ones and then multiplying first
two factors one obtains the required property of $\beta_U$. Hence the
bijective correspondence between $\beta_U$ and $\gamma$ is
established.

Finally, from Proposition~3.4, left-invariant $\omega_U$ is of the
canonical form with $i = (\eps\tens\id)\circ\omega_U$. Since $\omega_U(h) =
\Phi(Sh\o)\beta_U(\pi_\eps(h\t))\Phi(h\th) + \Phi(Sh\o)\extd_U\Phi(h\t)$
one easily finds that $i(h) = \gamma(h\o)\Phi(h\t)$, i.e. $i(h) =
\gamma(h\o)\tens 
h\t$, where $\gamma: H\to M$, $\gamma(h) =
(\eps\tens\id)\circ\beta(\pi_\eps(h))
+\eps(h)1_M$, for any $h\in H$.
\endproof

Therefore, for these $\beta_U$ we are in the setting of
Proposition~3.4 or Example~3.6 for the map $i$ constructed above.
The smallest horizontal right ideal in this case is
\begin{equation}
 \CQ_0=\span\{\gamma(q\o)q\t\la m\tens q\th h - \eps(m)\gamma(q\o
h\o)\tens q\t h\t \; | \; q\in \CQ, m\in M, h\in H\}.
\label{q0bic}
\end{equation}
 We see that
a choice of left-invariant $\omega_U$, $\CQ_{\rm
hor}\supseteq\CQ_0$ and left-covariant $\Omega^1(M)$ defines a
left-covariant $\Omega^1(P)$. The corresponding ideal is
$\CQ_P=\<i(\CQ)P,\CQ_MP\>$ where
\[ i(\CQ)P=\span\{\gamma(q\o)q\t\la m\tens q\th h|\ q\in\CQ, m\in M,
h\in H\}\supseteq\CQ_0.\]

\begin{example} Let $P=M\bicross H$ be viewed as a quantum principal
bundle. Let $\gamma$ obey the condition
in Proposition~6.1 and let $\Omega^1(M)$ be left $M$-covariant. Then
$P$ has a natural left-covariant
calculus $\Omega^1(P)$ such
that $\Omega^1_{\rm hor}=P(\extd M)P$ and
\[
\omega(h)=\Phi^{-1}(h\o)\beta(\pi_\eps(h\t))\Phi(h\th)
+\Phi^{-1}(h\o)\extd\Phi(h\t)\]
where $\beta:\ker\eps\to \Omega_P^1(M)$ is defined by
$\beta(h)=(S\gamma(h)\o)\extd\gamma(h)\t$.
\end{example}
\proof Since $\Omega^1(M)$ is assumed to be left-covariant, the
subbimodule $\CN_M$ generating $\Omega^1(M)$ is obtained from a right
ideal $\CQ_M\subset \ker\eps_M$. Since $M$ is a Hopf subalgebra of
$M\bicross H$, the left $M$-invariance of $\CN_M$ implies left
$P$-invariance of
$P\CN_M P$. The corresponding right ideal in $\ker\eps_P$ is  $\CQ_M
P$. Therefore we take $\CQ_{\rm hor}=\<\CQ_0,
\CQ_MP\>$ corresponding to $\CN_{\rm
hor}=\<\CN_0,P\CN_MP\>$ as in Example~3.7.  On the
other hand, we are also in the setting of Proposition~3.3 and take
$\omega_U$ in that strong form, as in Proposition~6.1. Note as in
Proposition~3.3. that the inherited differential structure
$\Omega^1_P(M)$ is not different from $\Omega^1(M)$ if $\CN_0\cap
\Omega^1M\subseteq \CN_M$.
\endproof

We now consider the simplest concrete setting of bicrossproducts,
 where $M=\C(\Sigma)$, $\Sigma$ a
finite set (as in Section~5) and $H=\C G$, $G$ a finite group. Here
$P$ is a bicrossproduct of the form $\C(\Sigma)\bicross
\C G$, regarded as a bundle. This is necessarily of the form
associated to a group factorisation
  $X=G\Sigma$. Then $\Sigma$ acts on $G$ and $G$ acts on $\Sigma$, by
$\la,\ra$ respectively, as defined by $sg=(s\la g)(s\ra g)$ in $X$.
The bicrossproduct $\C(\Sigma)\bicross \C G$ has the explicit form
\cmath{(\delta_s\tens g)(\delta_t\tens h)=\delta_{s\ra g,t}
(\delta_s\tens gh)\\
\Delta(\delta_s\tens g)=\sum_{ab=s}\delta_a\tens b\la g\tens
\delta_b\tens g,
\quad S(\delta_s\tens g)=\delta_{(s\ra g)^{-1}}\tens (s\la g)^{-1}}
for all $g,h\in G$ and $s,t\in \Sigma$. Note that the actions  $\la,\ra$ are
typically not effective. We define the subset
\[ Y=\{(g,s)|\ s\la g=g\}=\prod_{g\in G}I(g)\subseteq G\times \Sigma.\]
where $I(g)$ is the isotropy group of $g$. Here $Y$ necessarily
contains $\Sigma=I(e)$ as $(e,\Sigma)$ where $e\in G$ is the group
identity. From Proposition~5.1 we know that $\Omega^1(\C(\Sigma))$
correspond to $\Gamma\subset \Sigma\times\Sigma-{\rm diag}$. We
require this to be $\Sigma$-invariant. Finally, we know from
\cite{Ma:cla} that coirreducible bicovariant $\Omega^1(\C G)$
correspond to $(V,\lambda)$ where $V$ is an irreducible left
$G$-module and $\lambda\in P(V^*)$. The corresponding quantum
tangent space in \cite{Ma:cla} is spanned by $x_v=\lambda((\ )\la
v)-\lambda(v)1\in\C(G)$ with corresponding derivation $\del_{x_v}
g=x_v(g)g$ on group-like elements $g\in\C G$. Hence
\[ \CQ={\rm span}\{q\in\ker\eps|\
\eps\del_{x_v} q=0\}=\{q\in \ker\eps|\ \lambda(q\la v)=0\quad\forall v\in V\}\]
i.e. the kernel of the map $\ker\eps\to V^*$ provided by the action
$\la:\C G\tens V\to V$ composed with $\lambda$.

\begin{propos}  Left-invariant $\Omega^1(\C(\Sigma)\bicross\C G)$
are provided
by pairs $\gamma,S$, where $\gamma\in \C(Y)$ is a function such
that $\gamma(e,s)=1=\gamma(g,e)$ for all $s\in \Sigma, g\in G$, and
$S\subset \Sigma$,
$e\notin S$ is a subset.
The associated invariant connection is defined by 
\[ \beta_U(g)_{s,t}=\gamma(g,s^{-1}t)-1\]
where $\gamma$ is extended by zero to $X$. The minimal horizontal right ideal
is
\[ \CQ_0 = {\rm span}\{\sum_{q\in G}q_g\gamma(g)\delta_{s\ra g^{-1}}\tens gh|\
q\in\CQ, h\in G,\ e\ne s\in \Sigma\}
 +{\rm span}\{\sum_{g\in
G}q_g(\delta_e-\gamma(g))\tens g|q\in\CQ\}.\]
 If we take $\CQ_{\rm
hor}=\<\CQ_0,\C(S)\tens\C G\>$ then the resulting calculus has
\[ \CQ_P={\rm span}\{\sum_{q\in G}q_g\gamma(g)\delta_{s\ra g^{-1}}\tens gh|\
q\in\CQ, h\in G,\ e\ne s\in \Sigma\}+\delta_e\tens\CQ+\C(S)\tens\C G.\]
\end{propos}
\proof The coaction of $\C(\Sigma)$ on $\CG$ in the bicrossproduct is
$g\mapsto \sum_s s\la g\tens \delta_s$ (see \cite{Ma:book}). We
therefore require $\gamma:\C G\to
\C(\Sigma)$ i.e. $\gamma\in\C(G\times\Sigma)$ such that
\[\sum_s \gamma(g)\delta_s\tens s\la g=\gamma(g)\tens g\]
for all $g$. Evaluating at a fixed $s\in\Sigma$, this is
$\gamma(g,s)(s\la g-g)=0$ for all $s\in\Sigma$ and $g\in G$. This
gives the stated form of $\gamma$. Then
 \align{\equad&& \beta_U(g)=\sum_{s \in I(g)}
\gamma(g,s)\sum_{ab=s}\delta_{a^{-1}}\extd_U\delta_b
=\sum_{s\in I(g)}\gamma(g,s)\sum_{ab=s}\delta_{a^{-1}}
\tens\delta_b-\delta_{a^{-1}}\delta_b
\tens 1\\
&&=\sum_{s\in I(g)}\gamma(g,s)\sum_{ab=s}\delta_{a^{-1}}
\tens\delta_b-\gamma(g,e)1\tens 1}
which gives the formula for components of $\beta_U$ as stated.

Next, we require $\Omega^1(\C(\Sigma))$ to be left
$\C(\Sigma)$-invariant, i.e. that $\Delta_L \CN_M\subset M\tens
\CN_M$ where $M=\C(\Sigma)$ has coproduct $\Delta\delta_s
=\sum_{ab=s}\delta_a\tens\delta_b$.
As in Section~5 we take $\CN_M={\rm span}\{\delta_s\tens\delta_t|\
(s,t)\in\Gamma\}$. The invariance is then equivalent to $\Gamma$
stable under the diagonal action of the group $\Sigma$. Such
$\Gamma$ are of the form $\Gamma=\{(s,t)|\ s^{-1}t\in S\}$ for some
subset $S$ not containing the group identity $e$. The right ideal
$\CQ_M$ in this case is
\[ \CQ_M={\rm span}\{\delta_s|\ s\in S\}=\C(S).\]
{}From the form of the algebra structure of the bicrossproduct, it is
clear that $\CQ_MP=\C(S)\tens\C G$.

To compute $\CQ_0$ we consider elements in $\CQ$ of the form
$q=\sum_{g\in G} q_g g$ and the delta-function basis for
$\C(\Sigma)$ in the formula (\ref{q0bic}). Then
\align{&&\equad \CQ_0={\rm span}\{ \sum_{g\in G}q_g(\gamma(g)
g\la\delta_s-\delta_{s,e}\gamma(gh))\tens gh|\ q\in\CQ,\ h\in G,\
s\in \Sigma)\}\\ &&={\rm span}\{\sum_{g\in
G}q_g\gamma(g)\delta_{s\ra g^{-1}}\tens gh|\ q\in\CQ, h\in G,\ e\ne
s\in \Sigma\}\\&&+{\rm span}\{\sum_{g\in G}q_g(\delta_e-\gamma(g))\tens
g|q\in\CQ\}} where we consider the cases where $s=e$ and $s\ne e$
separately. In the former part we wrote $g\la\delta_s=\delta_{s\ra
g^{-1}}$ while in the case $s=e$ we change variables from
$qh=\sum_g q_ggh$ to $q$ since $qh\in \CQ$ for all $h$. We span
over $q\in \CQ$ after fixing $s\in M$ and $h\in G$.  The
computation of $\CQ_P$ is similar. \endproof

We demonstrate this construction now in some examples based on
finite cyclic groups. When $G=\Z_n=\<g\>$, for the representation
$V$ defining a calculus on $\C G$ we take the 1-dimensional
representation where the generator $g$ acts as $e^{2\pi\imath\over
n}$. Its character $\chi$  corresponds to a conjugacy class in
$\hat\Z_n$ if we take the view $\C \Z_n\isom\C(\hat\Z_n)$. The
corresponding quantum tangent space is spanned by
$x=\chi-1\in\C(\Z_n)$ with corresponding derivation $\del_x
g^a=x(g^a)g^a$ for $a\in\{0,\cdots,n-1\}$. Hence
\[ \CQ=\{q\in \C G|\
\eps(q)=0,\ \chi(q)=0\}=\{q_a=\sum_{b=0}^{n-1}
e^{2\pi\imath ab\over n} g^b|\ a=1,2,3,n-2\}\]
The remaining basis element $n^{-1}q_{n-1}$ of $\ker\eps$ is dual
to $x$ and can be identified with the unique normalised left-invariant
1-form in the calculus.

Likewise, for a calculus on $\C(\Sigma)$ where $\Sigma=\Z_m=\<s\>$
we take for left-invariant calculus the one defined by
$S=\{s^2,s^3,\cdots,s^{m-1}\}$. Since $\Sigma$ is Abelian,
left-invariant calculi are automatically bicovariant, and this is
the natural 1-dimensional bicovariant calculus
$\Omega^1(\C(\Sigma))$ associated to the generator $s\in\Sigma$.
The ideal $\CQ_M$ consists of all functions vanishing at $e,s$. The
element $\delta_s$ is dual to the quantum tangent space basis
element $s-e$ and can be identified with the unique normalised
left-invariant 1-form.

There are many factorisations of the form $\Z_n\Z_m$. We consider
one of the simplest, namely $S_3=\Z_2\Z_3$ (actually a semidirect
product) where $G=\Z_2=\<g\>$ and $\Sigma=\Z_3=\<s\>$. In terms of
permutations $\alpha,\beta$ obeying $\alpha^2=\beta^2=e$ and
$\alpha\beta\alpha=\beta\alpha\beta$, we write $g=\alpha$ and
$s=\alpha\beta$. The action $\la$ is trivial while $s\ra g=s^2$ and
$s^2\ra g=s$. The Hopf algebra $\C(Z_3)\lcross \C \Z_2$ is
6-dimensional with cross relations
\[ g\delta_e=\delta_eg,\quad g\delta_s=\delta_{s^2}g,\quad
g\delta_{s^2}=\delta_sg\]
and the tensor product coalgebra structure. The subset $Y$
is all of $G\times\Sigma$ and hence
\[ \gamma(e)=1,\quad \gamma(g)=\delta_e+\gamma_1\delta_s+\gamma_2\delta_{s^2}\]
for two parameters $\gamma_1,\gamma_2\in\C$.

\begin{example} For the cross product $P=\C(\Z_3)\lcross\C Z_2$ as above and
the
choice of the 1-dimensional $\Omega^1(\C \Z_2)$ and
$\Omega^1(\C(\Z_3))$ as above, we find $\Omega^1(P)$ is
3-dimensional corresponding to $\CQ_P={\rm
span}\{\delta_{s^2}\}\tens\C \Z_2$. It has basis of invariant forms
$\{\omega_0,\omega_1,\omega_2\}$ say and
\[\extd \delta_e=(\delta_{s^2}-\delta_e)\omega_1,
\quad\extd \delta_s=(\delta_e-\delta_s)\omega_1,\quad \extd g
=g(\omega_0-\omega_1+\omega_2)\]
and module structure
\[ \omega_0g=-g\omega_0,\quad \omega_1g=g\omega_2,\quad \omega_2g=g\omega_1\]
\[\omega_0\delta_{s^i}=\delta_{s^i}\omega_0,
\quad \omega_1\delta_{s^i}=\delta_{s^{i-1}}\omega_1,\quad \omega_2\delta_{s^i}
=\delta_{s^{i-1}}\omega_2.\]
The gauge field corresponding to $\gamma$ is
\[ \beta_U(g)=\pmatrix{0&\gamma_1&\gamma_2\cr
                       \gamma_2&0&\gamma_1\cr
                       \gamma_1&\gamma_2&0}\]
but the entries $\gamma_i$ do not affect the resulting calculus.
\end{example}
\proof The ideal $\CQ=0$ in this case, i.e. $\Omega^1(\C \Z_2)$ is being
taken here with the universal differential calculus, which is 1-dimensional
in the case of $\C\Z_2$. This is clear
from the point of view of a bicovariant calculus on $\C(\hat\Z_2)$.
Hence $\CQ_0=0$ as well, and we take $\CQ_{\rm hor}=\CQ_MP={\rm
span}\{\delta_{s^2}\}\tens\C \Z_2$ for all $\gamma$. According to
Proposition~3.5, $\CQ_P=\<\CQ_{\rm hor},i(\CQ)P\>=\CQ_MP$ as well
since $\CQ=0$. This gives the calculus $\Omega^1(P)$. It projects to
the  universal one
in the fibre direction and restricts to the initial calculus on the
base.

We now compute this 3-dimensional calculus explicitly. We recall
that $\Omega^1(P)=P\tens \ker\eps/\CQ_P$ as a left $P$-module by
multiplication by $P$, as a right $P$ module by $[h]u=u\o\tens [h
u\t]$ for $[h]\in\ker\eps/\CQ_P$ and $u\in P$. Here $[\  ]$ denotes
the canonical projection from $\ker\eps$. The exterior derivative
is $\extd u=u\o\tens u\t-u\tens 1$ projected to $\ker\eps/\CQ_P$.
In our case, a basis for the latter is
\[\omega_0=[\delta_e\tens(g-e)],\quad \omega_1=[\delta_s\tens e],\quad
\omega_2=[\delta_s\tens g].\]
Then $\extd (1\tens g)=1\tens g \tens [1\tens (g-e)]$ giving the
result as stated on identifying $g\equiv 1\tens g$ in $P$.
Moreover, $\extd (\delta_{e}\tens e)=\delta_e\tens e\tens
\delta_e\tens e+\delta_{s^2}\tens e\tens\delta_s\tens
e-\delta_e\tens e\tens 1\tens e=(\delta_{s^2}\tens e-\delta_e\tens
e)\tens [\delta_s\tens e]$ as stated, on identifying $\delta_s\tens
e\equiv\delta_s$ etc. Likewise, $\extd(\delta_s\tens
e)=\delta_s\tens e\tens\delta_e\tens e+\delta_{e}\tens
e\tens\delta_s\tens e-\delta_s\tens e\tens 1\tens e=(\delta_e\tens
e-\delta_s\tens e)[\delta_s\tens e]$ as stated.

Finally, we compute the right module structure as follows. For the
action on $\omega_0$ we have $[\delta_e\tens (g-e)](1\tens
g)=1\tens g \tens [(\delta_e\tens (g-e)) (1\tens g)]=-1\tens g\tens
[\delta_e\tens(g-e)]$ as stated. And $ [\delta_e\tens
(g-e)](\delta_{s^i}\tens e)=\sum_{a+b=i}
\delta_{s^a}\tens e\tens [(\delta_e\tens (g-e))(\delta_{s^b}\tens e)]
=\delta_{s^i}\tens e\tens [\delta_e\tens(g-e)]$ as stated. Only the
$b=0$ term in the sum contributes. For the action on $\omega_1$ we
have $[\delta_s\tens e](1\tens g)
=1\tens g\tens [(\delta_s\tens g)(1\tens g)]=[\delta_s\tens e]$.
And $[\delta_s\tens e](\delta_{s^i}\tens e)=\sum_{a+b=i}\delta_{s^a}
\tens e\tens [(\delta_s\tens e)(\delta_{s^b}\tens e)]=\delta_{s^{i-1}}
\tens e\tens [\delta_s\tens e]$ as only the $b=1$ term in the sum contributes.
Similarly for the action on $\omega_2$. As a left module the action
is free, i.e. we identify $(1\tens g)\tens\omega_0=g\omega_0$ etc.
\endproof

\begin{example} For the cross product $P=\C(\Z_3)\lcross\C Z_2$ as
above but the choice of zero differential calculus $\Omega^1(\C
\Z_2)$ and universal calculus $\Omega^1\C(\Z_3)$, we find $\Omega^1(P)$
is the zero calculus unless $\gamma_1\gamma_2=1$, when it is
2-dimensional. In the latter case, with  basis of invariant forms
$\{\omega_1,\omega_2\}$ we have
\[\extd \delta_e=(\delta_{s^2}-\delta_e)\omega_1
+\gamma_1(\delta_s-\delta_e)\omega_2,\quad\extd
\delta_s=(\delta_e-\delta_s)\omega_1
+\gamma_1(\delta_{s^2}-\delta_s)\omega_2,\quad \extd
g=(1-\gamma_1)g(\omega_2-\omega_1)\] and right module structure
\[  \omega_1g=g\omega_2,\quad \omega_2g=g\omega_1, \quad
\omega_1\delta_{s^i}=\delta_{s^{i-1}}\omega_1,
\quad \omega_2\delta_{s^i}=\delta_{s^{i-1}}\omega_2.\]
The restriction to $\Omega^1_P(\C(\Z_3))$ is a direct sum of the
 1-dimensional calculus associated to $s$ and the 1-dimensional
 calculus associated to $s^2$.
\end{example}
\proof If we take the zero differential calculus on
$\Omega^1(\C \Z_2)$, so $\CQ=\C(g-e)$, then
\[ \CQ_0={\rm span}\{\gamma_1\delta_s\tens gh-\delta_{s^2}\tens
h,\gamma_2\delta_{s^2}\tens gh-\delta_s\tens h|\ h\in\Z_2 \}.\]
Here the $s$ contribution to $\CQ_0$ is $\gamma(g)\delta_{s^2}\tens
gh-\gamma(e)\delta_s\tens h=\gamma_2\delta_{s^2}\tens
gh-\delta_s\tens h$. Similarly, the $s^2$ contribution is
$\gamma_1\delta_s\tens gh-\delta_{s^2}\tens h$. Finally, the $s^0$
contribution is $(\delta_e-\gamma(g))\tens g-(\delta_e-1)\tens
e=-\gamma_1\delta_s\tens g+\delta_{s^2}\tens
e-\gamma_2\delta_{s^2}\tens g+\delta_s\tens e$ is already
contained. This is 4-dimensional for generic $\gamma$ (in this case
$\CQ_0=\ker\eps\tens\C\Z_2$) but collapses to a 2-dimensional ideal
when $\gamma_1\gamma_2=1$.

If we take the universal calculus on $\C(\Z_3)$ so $\CQ_M=0$, we
have $\CQ_{\rm hor}=\CQ_0$ is 4-dimensional in the generic case or
2-dimensional in the degenerate case (note that if we took the
1-dimensional calculus on $\C(\Z_3)$ as before then $\CQ_{\rm
hor}=\ker\eps\tens\C\Z_2$ is 4-dimensional in either case).
Finally, $i(\CQ)=\gamma(g)\tens g-1\tens e=\delta_e\tens
(g-e)+\gamma_1\delta_s\tens g-\delta_{s^2}\tens
e+\gamma_2\delta_{s^2}\tens g-\delta_s\tens e$ so $\CQ_P=\<\CQ_{\rm
hor},i(\CQ)P\>=\ker\eps\tens\C\Z_2+\C\delta_e\tens(g-e)=\ker\eps$
is 5-dimensional except in the degenerate case when
$\gamma_1\gamma_2=1$. This means that the calculus on $P$ is the
zero one except in the degenerate case. In the degenerate case,
$\CQ_P={\rm span}\{q,qg,\delta_e\tens(g-e)\}$ is 3-dimensional,
where $q=\gamma_1\delta_s\tens g-\delta_{s^2}\tens e$ as a
shorthand.

In this degenerate case, a basis of $\ker\eps/\CQ_P$ is
\[ \omega_1=[\delta_s\tens e],\quad \omega_2=[\delta_s\tens g]\]
while in this quotient, $\delta_{s^2}\tens e=\gamma_1\delta_s\tens
g$ and $\delta_{s^2}\tens g=\gamma_1\delta_s\tens e$ instead of
zero as in the preceding example, while $\delta_e\tens(g-e)$ is now
zero in the quotient. The computations proceed as on the preceding
example with these changes, resulting in some extra terms with
$\gamma_1$ as stated. The restriction to $\C(\Z_3)$ has a part
spanned by $\omega_1$ which is the 1-dimensional calculus on
$\C(\Z_3)$ as in the preceding example and a part spanned by
$\omega_2$ which has a similar form when computed for
$\extd\delta_{s^2}$. \endproof

For a more complicated example one may take $S_3\times
S_3=\Z_6\dcross\Z_6$ in \cite{BGM:fin}, which is a genuine double
cross product with both $\la,\ra$ nontrivial.  Writing
$G=\Z_6=\<g\>$ and $\Sigma=\Z_6=\<s\>$, say, the actions of the
generators are by group inversion on the other group. Thus
\[ I(e)=I(g^3)=\Sigma,\quad I(g)=I(g^2)=I(g^4)=I(g^5)=\{e,s^2,s^4\}=\Z_3.\]
The space of allowed $\gamma$ is therefore 13-dimensional. The
bicrossproduct Hopf algebra $P=\C(\Z_6)\bicross\C Z_6$ in this case
is 36-dimensional. The results in this case are similar to the
situation above: for generic parameters one obtains the zero
calculus but for special values one obtains calculi on
$\Omega^1(P)$ restricting to non-universal calculi on the base.

Finally, one may apply Proposition~6.2 equally well in the setting
of Lie bicrossproducts. As shown in  \cite{Ma:mat} one has examples
$\C(G^{\star\rm op})\bicross U(\cg)$ for all simple Lie algebras
$\cg$. Here $G^{\star\rm op}$ is the solvable group in the Iwasawa
decomposition of the complexification of the compact Lie group $G$
with Lie algebra $\cg$. Such bicrossproduct quantum groups arise as
the actual algebra of observables of quantum systems, for example
the Lie bicrossproduct $\C(SU_2^{\star
\rm op})\bicross U(su_2)$ is the quantum algebra of observables of
a deformed top\cite{Ma:hop}\cite{Ma:book}. We consider this example
briefly. We take $\C(SU_2^{\star\rm op})$ as described by
coordinates $\{X_i\}$ and $(X_3+1)^{-1}$ adjoined, and a usual
basis $\{e_i\}$ of $su_2$. Then the bicrossproduct is (see
\cite{Ma:book})
\[ [X_i,X_j]=0, \ \  \Delta X_i=X_i\tens 1+(X_3+1)\tens X_i,\ \
\eps X_i=0,\ \   SX_i=-{X_i\over X_3+1}.\]
\cmath{[e_i,e_j]=\eps_{ijk}e_k,\quad [e_i,X_j]
=\eps_{ijk}X_k-\h\eps_{ij3}{X^2\over X_3+1},\quad \eps e_i=0,\\
\Delta e_i=e_i\tens {1\over X_3+1}+e_3\tens {X_i\over X_3+1}+1\tens e_i,\quad
Se_i
=e_3X_i-e_i(X_3+1).}
For a differential calculus $\Omega^1(\C(SU_2^{\star\rm op}))$ we
have a range of choices including the standard commutative one.
Others are ones with quantum tangent space given by jet
bundles\cite{Ma:cla}. For $\Omega^1(U(su_2))$ one may follow a
similar prescription to $\Omega^1(\C G)$: if $V$ is an irreducible
representation and $\lambda\in P(V^*)$ then $\CQ=\{q\in U(su_2)|\
\eps(q)=0, \ \lambda(q\la v)=0,\ \forall v\in V\}$. A natural choice
is $V$ a highest weight representation and $\lambda$ the conjugate
to the highest weight vector. Finally, we consider the possible
$\gamma$. Note first of all that in a von-Neumann algebra setting
one may consider group elements $g\in SU_2$ much as in
Proposition~6.3. From the explicit formulae for the action of
$su_2^{\star\rm op}$ on $SU_2$ in \cite{Ma:mat}\cite{Ma:book}, one
then sees that at least near the group identity,
\[ I(g)=\{\exp t (f_3-{\rm Rot}_g(f_3))|\ t\in\R\}\]
where $\{f_i\}$ are the associated basis of the Lie algebra
$su_2^{\star\rm op}$ and ${\rm Rot}$ is the action of $SU_2$ by
rotations of $\R^3$. Hence $\gamma$ should be some form of
distribution on $SU_2\times SU_2^{\star\rm op}$ such that
$\gamma(g)$ has support in the line $I(g)$. This suggests that in
our algebraic setting one should be able to construct a variety of
$\gamma:U(su_2)\to
\C(SU_2^{\star\rm op})$ order by order in a basis of $U(su_2)$. Thus,
at the lowest order the coaction of $\C(SU_2^{\star\rm op})$
is\cite{Ma:book}
\[ e_i\bo\tens e_i\bt=e_i\tens (X_3+1)^{-1}+e_3\tens X_i(X_3+1)^{-1}\]
and the condition for $\gamma$ in Proposition~6.2 becomes
\[\gamma(e_i)X_3-\gamma(e_3)X_i=0.\]
This has solutions of the form $\gamma(e_i)=f_i(X)X_i$ for any
functions $f_i\in\C(SU_2^{\rm op})$. After fixing $\gamma$ and the
base and fibre differential calculi one may obtain left-invariant
differential calculus on $P=\C(SU_2^{\star
\rm op})\bicross U(su_2)$ and a connection on it as a quantum principal bundle.
The detailed analysis will be considered elsewhere.

We note that as a semidirect product one could also think of this
bicrossproduct as a deformation of the 3-dimensional Euclidean
group of motions (one may introduce a scaling of the $X_i$ to
achieve this). In the 3+1 dimensional version of this same
construction one has the $\kappa$-deformed Poincar\'e algebra as
such a bicrossproduct\cite{MaRue:bic}. Proposition~6.2 therefore
provides in principle a general construction for left-invariant
calculi on these as well. At the moment, only some examples are
known by hand \cite{KMS:bic}.  Moreover,  affine quantum groups
such as $U_q(\hat{su}_2)$ may be considered as cocycle
bicrossproducts $\C[c,c^{-1}]\bicross U_q(Lsu_2)$ where
$U_q(Lsu_2)$ is the level zero affine quantum group (quantum loop
group) and $c$ is the central charge generator, see \cite{Ma:aff}.
The quantum Weyl groups provide still more examples of cocycle
bicrossproducts\cite{MaSoi:wey}. All of these and their duals may
be treated as (trivial) quantum principal bundles by similar
methods to those above.

\subsection{Biproducts, bosonisations and the quantum double}

Let $H$ be a Hopf algebra with (for convenience) bijective
antipode. A braided group in the category of crossed modules means
$B$ which is an algebra, a coalgebra and a crossed $H$-module (i.e.
a left $H$-module and left $H$-comodule in a compatible way) with
all structure maps intertwining the action and coaction of $H$ and
with the coproduct $\und\Delta:B\to B\und\tens B$ a homomorphism in
the braided tensor product algebra structure $B\und\tens B$. This
is basically the same thing as a braided group in the category of
$D(H)$-modules where $D(H)$ is Drinfeld double in the
finite-dimensional case. One knows from the braided
setting\cite{Ma:skl} of \cite{Rad:str} that every such braided
group has an associated Hopf algebra $B\lbiprod H$ as cross product
and cross coproduct. Moreover, $\pi(b\tens h)=\und\eps(b)h$ defines
a projection $B\lbiprod H\to H$ split by Hopf algebra map
$j(h)=1\tens h$. All split Hopf algebra projections to $H$ are of
this from.

We can clearly view such $B\lbiprod H$ as principal bundles of the
homogeneous type, with $j$ defining a canonical
connection\cite{BrzMa:gau}. The homogeneous bundle coaction is
$\Delta_R(b\tens h)=b\tens h\o\tens h\t$ so that $M=B$. Since $j$
is a coalgebra map we can also take $\Phi=j$ as a trivialisation,
i.e. the bundle is trivial. Both Propositions~3.3 and~3.5 apply in
this case. From the former, we know that any
$\beta_U:\ker\eps\to\Omega^1B$ and any $\Omega^1(H)$ yields a
calculus $\Omega^1(B\lbiprod H)$. We now use Proposition~3.5 to
study which of these are left-invariant.

Note that $B$ as a braided group coacts on itself via the braided
coproduct. This is the braided left regular coaction. This extends
to $B\und\tens B$ as a braided tensor product coaction, via the
braiding $\Psi(v\tens w)=v\bo\la w\tens v\bt$ of the category of
crossed modules. So
\[ \und\Delta_L(b\tens c)=b\Bo\Psi(b\Bt\tens c\Bo)\tens c\Bt
=b\Bo(b\Bt\bo\la c\Bo)\tens b\Bt\bt\tens c\Bt\]
and this restricts to a left $B$-coaction on $\Omega^1(B)$. The
calculus $\Omega^1(B)$ is braided-left covariant if its associated
ideal $\CN_B$ is stable under $\und\Delta_L$.

\begin{propos} Strong left-invariant connections $\omega_U$ on
$B\lbiprod H(B,H,j)$ are in
1-1 correspondence with the maps $\beta_U:\ker\eps_H\to\Omega^1B$
which are left $B$-invariant
(under the braided coproduct) and intertwine the left $H$-coaction on
$\Omega^1B$ with the left-adjoint coaction of $H$, i.e.
$$
\Delta_H(\beta_U(h)) = h\o Sh\th \tens \beta_U(h\t).
$$
Moreover, the $\beta_U$ are in 1-1 correspondence
with $\gamma:\ker\eps_H\to \ker\und\eps$ which are intertwiners of the left
adjoint coaction of $H$, i.e.
\[\Delta_H(\gamma(h)) = h\o Sh\th \tens \gamma(h\t). \]
The
correspondence is via
\[ \beta_U(h) = \und S\gamma(h)\Bo\extd_U\gamma(h)\Bt .\]
The corresponding $\omega_U$ is the  canonical connection for the splitting
$i(h)= \gamma(\pi_\eps(h\o))\tens h\t + 1\tens h$
\label{prop.biprod}
\end{propos}
\proof
$B\lbiprod H(B,H,j)$ is a trivial bundle with trivialisation
$j$. Given strong connection  $\omega_U$ one associates to
it the unique map
$\beta_U :\ker\eps_H\to \Omega^1 B$ given by
$\beta_U(h) = j(h\o)\omega_U(\pi_\eps(h\t))j(Sh\th)) + j(h\o)\extd_Uj(Sh\t)$. Since $j$ is a Hopf algebra map  the left coaction
$\Delta_L$ of
$B\lbiprod H$ on $\beta_U(h)\in\Omega^1 B\in \Omega^1  B\lbiprod H $
can be easily computed using the fact the $\omega_U$ is left-invariant
\begin{eqnarray*}
\Delta_L(\beta_U(h)) & = &
j(h\o)j(Sh\fiv)\tens j(h\t)\omega_U(\pi_\eps(h\th))j(Sh\fo) +
j(h\o)j(Sh\fo)\tens j(h\t)\extd_U j(Sh\th)\\
& = & j(h\o Sh\th)\tens \beta_U(h\t) = 1\tens h\o Sh\th\tens
\beta_U(h\t),
\end{eqnarray*}
where we also used the fact that $\ker\eps_H$ is invariant under the
left adjoint coaction. On the other hand  the left coaction of
$B\lbiprod H$ on $B\tens B\subset (B\lbiprod
H)^{\tens 2}$ is
\align{&&\equad \Delta_L(b\tens c)
=(b\Bo\tens b\Bt\bo)(c\Bo\tens c\Bt\bo)\tens b\Bt\bt\tens c\Bt\bt\\
&&=b\Bo (b\Bt\bo\o\la c\Bo)\tens b\Bt\bo\t c\Bt\bo\tens b\Bt\bt\tens
c\Bt\bt\\ &&= b\Bo (b\Bt\bo\la c\Bo)\tens b\Bt\bt\bo c\Bt\bo\tens
b\Bt\bt\bt\tens c\Bt\bt=(\id\tens\Delta_H)\und\Delta_L(b\tens c)}
where $\Delta_H$ is the left tensor product coaction of $H$ on $B\tens
B$. We used the comodule property for the $H$-coaction for the
third equality. Therefore we have just found that
\begin{equation}
(\id\tens\Delta_H)\und\Delta_L \beta_U(h) = 1\tens h\o Sh\th\tens
\beta_U(h\t).
\label{adj}
\end{equation}
Applying $\id_B\tens\eps_H\tens\id_B\tens\id_B$ to both sides of
(\ref{adj}) we find $\und\Delta_L \beta_U(h) = 1\tens \beta_U(h)$,
i.e. $\beta_U$ is left-invariant with respect to the (braided) left
coaction of $B$. Using this left-invariance we can compute (\ref{adj})
further to find
\[1\tens h\o Sh\th\tens
\beta_U(h\t) = (\id\tens\Delta_H)\und\Delta_L \beta_U(h) =
(\id\tens\Delta_H)(1\tens \beta_U(h)) = 1\tens \Delta_H(\beta_U(h)),\]
which is the required intertwiner property of $\beta_U$.

Conversely, assume that $\beta_U:\ker\eps_H\to \Omega^1 B$ is left
$B$-invariant and an intertwiner for the left adjoint coaction of $H$.
One then immediately finds for any $h\in\ker\eps_H$
\[\Delta_L(\beta_U(h)) = (\id\tens\Delta_H)\und\Delta_L \beta_U(h) =
1\tens h\o Sh\th\tens \beta_U(h\t).\]
Using this fact one computes
\begin{eqnarray*}
\Delta_L(\omega_U(h)) &=&
\Delta_L(j(Sh\o)\beta_U(\pi_\eps(h\t))j(h\th) + Sj(h\o)\extd_Uj(h\t))
\\
& = & 1_B\tens Sh\t h\th Sh\fiv h\six \tens
j(Sh\o)\beta_U(\pi_\eps(h\fo))j(h\sev) +1_{B\lbiprod H}\tens
Sj(h\o)\extd_Uj(h\t)\\
& = &  1_{B\lbiprod H}\tens
j(Sh\o)\beta_U(\pi_\eps(h\t))j(h\th)+1_{B\lbiprod H}\tens
Sj(h\o)\extd_Uj(h\t)\\
& = & 1_{B\lbiprod H}\tens\omega_U(h),
\end{eqnarray*}
so the connection corresponding to $\beta_U$ is left $B\lbiprod
H$-invariant.

Similar arguments as in the proof of Proposition~3.5 show that any
left $B$-invariant $\beta_U:\ker\eps_H\to \Omega^1B$ can be expressed
in the form
$\beta_U(h) = \und S\gamma(h)\Bo\extd_U\gamma(h)\Bt$, where $\gamma:
\ker\eps_H\to \ker\und\eps$ is given by $\gamma =
(\und\eps\tens\id)\circ\beta_U$. Since $\beta_U$ is an intertwiner for
the left adjoint coaction of $H$ we have
\[ \gamma(h)\Bo\bo \gamma (h)\Bt\bo \tens \und
S\gamma(h)\Bo\bt\extd_U\gamma(h)\Bt\bt = h\o Sh\th\tens
\und S\gamma(h\t)\Bo\extd_U\gamma(h\t)\Bt ,\]
where we used that $\und S$ is a left $H$-comodule map. Due to the
form of $\extd_U$, the above equality is in $H\tens B\tens B$. Applying
$\und\eps$ to the middle factor and using
the fact that $B$ is an $H$-comodule coalgebra one finds
\[\Delta_H(\gamma(h)) = h\o Sh\th \tens \gamma(h\t), \]
i.e. the required intertwiner property. Conversely, given
$\gamma:\ker\eps_H\to \ker\und\eps$ which is an intertwiner for the
left adjoint coaction we find
\begin{eqnarray*}
\Delta_H(\beta_U(h)) & = & \Delta_H(\und
S\gamma(h)\Bo\extd_U\gamma(h)\Bt) \\
& = & \gamma(h)\Bo\bo \gamma (h)\Bt\bo \tens \und
S\gamma(h)\Bo\bt\extd_U\gamma(h)\Bt\bt\\
& = & \gamma(h)\bo \tens \und
S\gamma(h)\bt\Bo\extd_U\gamma(h)\bt\Bt\\
& = & h\o Sh\th\tens \und S\gamma(h\t)\Bo\extd_U\gamma(h\t)\Bt = h\o
Sh\th\tens \beta_U(h\t),
\end{eqnarray*}
as required. Finally, if the map $\beta_U$ is expressed in terms of the
map $\gamma$ then the canonical splitting $i$ corresponding to
$\omega_U$ and given by $i = (\eps\tens\id)\circ\omega_U$ comes out as
stated in the proposition. \endproof

Notice that the map $\gamma :\ker\eps_H\to\ker\und\eps$ defined in
Proposition~6.6 can be  uniquely extended
to the map $\gamma : H\to B$ by requiring $\gamma(1) =1$. Then
$\underline{\eps}\circ\gamma = \eps$ and $i =
(\gamma\tens\id)\circ\Delta$. Therefore for these $\beta_U$ we are in
the setting of
Proposition~3.5 or Example~3.7 for the stated map $i$. The smallest
horizontal right ideal is
\[ \CQ_0=\span\{\gamma(q\o)q\t\la b\tens q\th h - \und\eps(b)\gamma(q\o
h\o)\tens q\t h\t \; | \; q\in \CQ, b\in B, h\in H\}.\]
  We see
that a choice  of left-invariant $\omega_U$, $\CQ_{\rm hor}$
and left-covariant $\Omega^1(B)$ yields a suitably left-invariant
$\Omega^1(B\lbiprod H)$.

\begin{example} Let $P=B\lbiprod H$ be viewed as a quantum principle bundle
as above with trivialisation $j(h) = 1\tens h$. Let $\gamma$ obey the
condition in Proposition~6.6 and
let $\Omega^1(B)$ be braided left $B$-covariant and $H$-covariant.
Then $P$ has a natural left-covariant calculus $\Omega^1(B\lbiprod
H)$ such that $\Omega^1_{\rm hor}=P(\extd B)P$ and
\[ \omega(h) = j(Sh\o)\beta(\pi_\eps(h\t)) j(h\th) + j(Sh\o)\extd j(h\t)\]
for $\beta:\ker\eps\to\Omega^1_P(B)$ defined by $\beta(h)=\und
S\gamma(h)\Bo\extd \gamma(h)\Bt$ is a
connection on it. Here $\Omega^1_P(B) = \pi_{\CN}(\Omega^1B)$, where
$\pi_{\CN} : \Omega^1P\to \Omega^1(P)$ is the canonical surjection.
\label{example.biprod}
\end{example}
\proof We take $\CQ_{\rm hor}=\<\CQ_0,\CQ_BP\>$ corresponding to $\CN_{\rm hor}
=\<\CN_0,P\CN_BP\>$ as in Example~3.6.
 Since $\Delta_L(b\tens c)=(\id\tens\Delta_H)\und\Delta_L(b\tens c)$
for any $b\tens c\in B\tens B$ viewed inside $(B\lbiprod
H)^{\tens 2}$ on the left hand side, it is clear that if $\Omega^1(B)$ is
defined by $\CN_B$ which is both $\und\Delta_L$ and $\Delta_H$
covariant then it is covariant under the $\Delta_L$ coaction of
$B\lbiprod H$. We then use the preceding Proposition~6.6 to establish
that $\omega_U(h) = j(Sh\o)\beta_U(\pi_\eps(h\t))j(h\th) + j(Sh\o)\extd_Uj(h\t)$, where $\beta_U(h) = \und
S\gamma(h)\Bo\extd_U \gamma(h)\Bt$ is a left-invariant
connection. We extend $\gamma$ to the whole of $H$ by setting
$\gamma(1) =1$. Then the  corresponding splitting is $i(h) =
\gamma(h\o)\tens h\t$ and we  construct a
left-covariant calculus $\Omega^1(B\lbiprod
H)$ by taking $\CQ_{B\lbiprod
H} = \<\CQ_{\rm hor}, i(\CQ)(B\lbiprod
H)\> = {\rm span}\{q_Bb\tens h, \gamma(q\o)q\t\la b\tens q\th h \; | \;
q\in \CQ, q_B\in \CQ_B, b\in B, h\in H\} $,
 as in Proposition~\ref{prop.homog}. \endproof

Bosonisation may be viewed as a special kind of biproduct, albeit
originating\cite{Ma:bos} from other considerations than
\cite{Rad:str}. Here $H$ is a dual quasitriangular Hopf
algebra\cite{Ma:bg}cf\cite{Dri} and $B$ a braided group in its
category of (say) left comodules. It has a bosonisation $B\lbiprod
H$ where the required action is induced by evaluating against the
dual quasitriangular structure, see \cite{Ma:war95} and
cf\cite{BrzMa:gau} (where the example of the quantum double as a
bundle was emphasised).  Explicitly,
\cmath{ (b\tens h)(c\tens g)=bc\bt\tens h\t g \CR(c\bo\tens h\o), \quad
\Delta(b\tens h)=b\Bo\tens b\Bt\bo h\o\tens b\Bt\bt\tens h\t}
where $\Delta_Lb=b\bo\tens b\bt$ is the coaction of $H$ (summation
understood).

To give a concrete application, we take $B=\underline{H}$ to be a
braided version of $H$ obtained by transmutation\cite{Ma:bg}.
$\underline{H}$ is then a left $H$-module coalgebra with the
coaction provided by left adjoint coaction, i.e. $\Delta_H (h) =
h\o Sh\th\tens h\t$. It is clear that the map $\gamma :
\underline{H}\to H$, $\gamma(h) = h$ satisfies all the requirements
of Proposition~\ref{prop.biprod} therefore we have
\begin{propos}
Let $B = \underline{H}$, $P=\underline{H}\lbiprod H$ and $\gamma =
\id$. Let the braided left-covariant calculus on $\underline{H}$ be
generated by $\CQ_{\underline{H}}\subset\ker\underline{\eps}$. Assume that the
ideal $\CQ\subset
\ker\eps_H$ is generated by $\{q_i\}_{i\in I}$. Then the corresponding
right ideal $\CQ_P\subset\ker\eps_P$ defining left-covariant
calculus on $P=\underline{H}\lbiprod H$ as in
Example~\ref{example.biprod} is generated by $\{\Delta q_i\}_{i\in
I}$ and the generators of $\CQ_{\underline{H}}$. The induced
calculus $\Omega^1_P(\und{H})$ is generated by $\<\CQ,
\CQ_{\underline{H}}\>$.
\label{prop.biprod.con}
\end{propos}
\proof The braided product ${\underline\cdot}$ in $\underline{H}$ is related to
the original product in $H$ by
\[ gh = g\o{\underline\cdot}h\t \CR(h\o Sh\th \tens g\t).\]
Since $\gamma$ is an identity map one finds the corresponding
splitting $i =\Delta$. For any $g,h\in\ker\eps$ we find
\begin{eqnarray*}
i(g)i(h) & = & (g\o\tens g\t)(h\o\tens h\t) = g\o \underline{\cdot}
h\t\CR(h\o Sh\th \tens g\t)\tens g\th h\t \\
& = & g\o h\o\tens g\t h\t = i(gh).
\end{eqnarray*}
This implies that if $\CQ$ is generated by $\{q_i\}_{i\in I}$ as a
right ideal in $H$ then $\CQ_P = \<\CQ_{\underline{H}}P, i(\CQ)P\>$ is
generated by generators of $\CQ_{\underline{H}}$ and
$\{i(q_i)\}_{i\in I}$ as a right ideal in $P$. Since $i$ is the same
as the 
coproduct in $H$, the assertion follows.

To derive the induced calculus on $\underline{H}$ first note that $\ker\eps_P =
\ker\underline{\eps}\tens 1 \oplus \underline{H}\tens \ker\eps_H$, where the
splitting
is given by  the projection $\Pi :\ker\eps_P \to \ker\underline{\eps}\tens
1$, $\Pi = \pi_\eps\tens\eps$. The differential structure on $\underline{H}$ is
determined by the image of $\CQ_P$ under this projection. Clearly
$\Pi(\CQ_{\underline{H}}P) =\CQ_{\underline{H}}$. Furthermore, for any
$b\in \underline{H}, h\in H, q\in\CQ$ we find
\begin{eqnarray*}
\Pi((q\o\tens q\t)(b\tens h)) & = & \Pi(q\o\underline{\cdot}b\t\CR(b\o
Sb\th \tens q\t) \tens q\th h) \\
& = &\Pi(q\o b\tens q\t h) =
 \pi_\eps(qb)\eps(h) = qb\eps(h).
\end{eqnarray*}
Therefore $\CQ_P$ restricted to $\underline{H}$ coincides with
 $\<\CQ_{\underline{H}},\CQ\>$ as stated.
 \endproof

In particular, if $\CQ_{\underline{H}}$ in the preceding
proposition is chosen to be trivial and thus the calculus on
$\underline{H}$ to be the universal one, the induced calculus is
non-trivial and is a braided version of  the calculus on $H$.
Notice also that since\cite{Ma:book}
\[ \underline{H}\lbiprod H \cong H\dcross H \cong D(H)^*,\]
(the latter in the case where $\CR$ is factorisable) the above
construction gives a natural left-covariant differential structure
on the double cross products $H\dcross H$ (see\cite{Ma:poi}) and
the duals of the Drinfeld double $D(H)$ (see \cite{Dri}). For
example, if $H=A(R)$, is a matrix quantum group corresponding to a
regular solution of the quantum Yang-Baxter equation (suitably
combined with q-determinant or other relations) then $\underline{H}
= B_L(R)$, the left handed version of the corresponding matrix braided
group\cite{Ma:exa}. The latter is generated by the matrix $\bf u$
subject to the left-handed braided matrix relations
$$
R{\bf u}_1R_{21}{\bf u}_2 ={\bf u}_2R{\bf u}_1R_{21}
$$
and suitable braided determinant or other relations. If $\bf t$
denotes the matrix of generators of $A(R)$ then the cross relations
in $B_L(R)\lbiprod A(R)$ are given by
$$
{\bf t}_1{\bf u}_2 = R_{21}{\bf u}_2R_{21}^{-1} {\bf t}_1
$$
as the left-handed version of the formulae in \cite{Ma:mec}. The isomorphism
with $A(R)\dcross A(R)$ as generated by $\bf s$ and $\bf t$ say (and the
cross relations $R{\bf s}_1{\bf t}_2={\bf t}_2{\bf s}_1R$) is ${\bf s}
={\bf u}{\bf t}$ and the $\bf t$ generators identified, as the appropriate
left-handed version of \cite{Ma:poi}.  In particular, taking $R$ to be the
standard $SU_q(2)$ R-matrix and
${\bf u} = \pmatrix{a& b \cr c& d}$, the relations for
 $BSU_q(2)=\underline{SU_q(2)}$ come out as
\[ da =ad, \quad cd = q^2 dc, \quad db = q^2bd, \quad bc = cb
 +(q^{-2}-1)(a-d)d,\]
\[ ac = ca + q^{-2}(1-q^{-2})cd, \quad ab = ba + (q^{-2}-1)bd, \quad ad-q^2bc=1
\]
The corresponding bosonisation $BSU_q(2)\lbiprod SU_q(2)$ is
isomorphic to the quantum Lorentz group $SU_q(2)\dcross SU_q(2)$
and thus Proposition~\ref{prop.biprod.con} allows one to construct
a differential calculus on the quantum Lorentz group. Moreover, as
explained in \cite{Ma:mec}, the bosonisation form of the quantum
double is  quite natural if we would like to regard it as a
q-deformed quantum mechanical algebra of observables or `quantum
phase space'. It is therefore natural to build its differential
calculus from this point of view.

Finally, one has braided covector spaces $V^*(R',R)$ in the
category of left $\tilde{A(R)}$-comodules, with additive braided
group structure. Here $\tilde{A(R)}$ denotes a dilatonic extension.
The bosonisation $V^*(R',R)\lbiprod \tilde{A(R)}$ has been
introduced in \cite{Ma:poi} as a general construction for
inhomogeneous quantum groups such as the dilaton-extended
$q$-Poincar\'e group $\R_q^n\lbiprod\tilde{SO_{q}(n)}$. The
detailed construction of the required intertwiner map $\gamma$ in
this case will be addressed elsewhere. We note only that in the
classical $q=1$ case it can be provided by a map $\gamma:\R^n\to
SO(n)$ such that $\gamma(g.x)=g\gamma(x)g^{-1}$ for all $g\in
SO_n$. For example, for $n=3$ the map $\gamma(x)=\exp(x)$ has this
property, where $x$ is viewed in $so_3$ by the Pauli matrix basis
and exponentiated in $SO(3)$. The $q$-deformed version of such maps
should then allow the application of the above methods to obtain
natural left-invariant calculi on inhomogeneous quantum groups as
well.


\end{document}